# The 42 Assessors and the Box-Kites they fly:

## Diagonal Axis-Pair Systems of Zero-Divisors in the Sedenions' 16 Dimensions

Robert P. C. de Marrais

**Abstract**. G. Moreno's abstract depiction of the Sedenions' normed zero-divisors as homomorphic to the exceptional Lie group G2, is fleshed out by exploring further structures the A-D-E approach of Lie algebraic taxonomy keeps hidden. A break-down of table equivalence among the half a trillion multiplication schemes the Sedenions allow is found; the 168 elements of PSL(2,7), defining the finite projective triangle on which the Octonions' 480 equivalent multiplication tables are frequently deployed, are shown to give the exact count of primitive unit zero-divisors in the Sedenions. (Composite zero-divisors, comprising all points of certain hyperplanes of up to 4 dimensions, are also determined.) The 168 are arranged in point-set quartets along the 42 Assessors (pairs of diagonals in planes spanned by pure imaginaries, each of which zero-divides only one such diagonal of any partner Assessor). These quartets are multiplicatively organized in systems of mutually zero-dividing trios of Assessors, a D4-suggestive 28 in number, obeying the 6-cycle crossover logic of trefoils or triple zigzags. 3 trefoils and 1 zigzag determine an octahedral vertex structure we call a box-kite -- seven of which serve to partition Sedenion space. By sequential execution of proof-driven production rules, a complete interconnected box-kite system, or Seinfeld production (German for field of being; American for 1990's television's Show About Nothing), can be unfolded from an arbitrary Octonion and any (save for two) of the Sedenions. Indications for extending the results to higher dimensions and different dynamic contexts are given in the final pages.

If, upon reading this paper, you wish to be included in a mailing list I'm contemplating on the subject, then send e-mail To: rdemarrais@alum.mit.edu and include "[42]" in the subject line.

**{ 0 } Introduction:** Some two centuries ago, de Moivre, Argand, Gauss, et al., determined that square roots of negative unity could be domesticated (and rendered geometric) by the introduction of cyclical 2-D numbers, ironically (given their futurity) deemed "imaginary." Since then, a succession of "period doublings," to borrow the Chaos theorists' term, has produced ever more exotic numbers – but at the cost of a loss of familiar behavior with each escalation. Already with the Imaginaries, the angular nature of the exponent made a product's "winding number" indeterminate. With the doubling to 4-D effected by Hamilton's Quaternions, three imaginary axes now could serve to map the composition of forces in the 3-D space of mechanics, with the fourth axis of Reals relating to time or tension . . . but at the cost of commutation. (Also lacking, of course, in the physical world such "vector mechanics" would model). The redoubling to the 8-D Octonions of Cayley-Dickson algebra required surrendering associativity – and swallowing the cause-effect indeterminacy such "parenthesis-placement uncertainty" implied. Hurwitz' century-old proof, meanwhile, showed that further twofold "blow-ups" in dimension of the number system would require the breakdown of the very notion of a division algebra: divisors of zero, for starters, would be let in at the same door the distributive property would exit through.

The development of quantum mechanics has made it clear that the benefits of such redoublings have, on the whole, outweighed the costs. The semi-simple Lie-group taxonomizing associated with Dynkin's (and, among geometers, Coxeter's) work has meanwhile allowed these developments to be tamed in a benign set of ball-and-stick diagrams, associated with neutral-sounding letters and indices. The complex numbers are A2, having the structure of 2-D closest packing patterns and the hexagonal arrays of Eisenstein integers generated by the cube roots of the positive and negative units. The 24 integer Quaternions conform to the "Feynman checkerboard" pattern of D4; and – most spectacularly of all – the closest packing of spheres in 8 dimensions, and the exceptional E8 related to icosahedral symmetry, correspond precisely to the 240-element integer arithmetic Coxeter elicited for Octonions.

Guillermo Moreno's recent (1997) discovery of the $G_2$ homomorphism governing the zero divisors in the 16-D Sedenions[1] would seem, on the surface, to be the latest triumph of this "A-D-E" approach. And indeed, having conquered virtually all taxonomic work in all areas of mathematics where finite typologies are pertinent, this is what we've come to expect. It's been only a quarter century since the great Russian dynamicist Vladimir Arnol'd discovered that "A-D-E" could house much more than the "particle zoos" of the physicists (and the complete catalog of Coxeter's kaleidoscopes). The classification of singularities begun by René Thom's Catastrophe Theory could be engulfed and extended by it, an embarrassment of riches leading Arnol'd to speak of the "A-D-E *Problem*."[2] Today, there is likely no item in the Coxeter-Dynkin "alphabet soup" that lacks a surprising realization as an unfoldable singularity with boundaries or other constraints (or a projective-pencil "alcove" in Tits' theory of buildings, or a class of quivers, light caustics, shock-waves or self-correcting error codes). Moreno's $G_2$ would thereby seem to be another indicator of what Arnol'd called the A-D-E-induced "mysterious unity of all things."[3]

But Moreno discovered a *homomorph*ism -- a "blow-up" of an exact correspondence – and the "blow-ups" in the history of number theory have all entailed the loss of something. We know now that among the numerous copies of the Octonions floating in the extra *lebensraum* of the Sedenion's 16-D, we can find divisors of zero harbored among them according to $G_2$'s 12-fold (or 14-fold, depending upon the questions asked) logic. But how many such distinct sets of $G_2$'s are there? Of what do these sets in fact consist? What is the *dynamic* meaning of zero-division in their context?

**{ I } Names, Notations, Conventions, Initial Constructions:**

As Moreno makes clear, a Sedenion zero divisor (unlike those in standard low-dimensioned Clifford algebras, typically generated by mutually annihilating idempotents with real components) is composed of purely imaginary units. "Doubly pure," in his terminology, as he implicitly concerns himself only with what will be called herein "primitive" zero-divisors, comprised of diagonals in planes spanned by ***pairs*** of imaginary units. (We will see that additive sums of up to four such primitives also exist, and for such composites it will no longer make much sense to focus, as Moreno does, on "normed unit" divisors: hyperplanes of up to 4 dimensions will be found, all of whose points are zero-divisors, so that the infinity of points on their unit spheres would all qualify for "normed unit" status.)

The Sedenions' primitive zero-divisors will be referred to by a special name, "Assessors," and the systems of these which mutually zero-divide each other (unlike the mutually annihilating ***pairs*** familiar from Clifford algebra[4]) will be called "***Co-Assessor trios***." And while technically there are four "normed unit" zero-divisors on the extended X that is an Assessor, in fact plus or minus differences along one or the other diagonal have no effect on the basic dynamics. I prefer, therefore, to think of each Assessor's "X" as made of two infinite ($S^0$-normed) lines, one an orthographic slash, the other a backslash, all points of which are zero-divisors with respect to one or the other such line of a Co-Assessor: ‘/’ + ‘\’ = ‘X’. In the spirit of the bit-twiddling to come, I call these paths “***UNIX*** ” and “***DOS*** ” respectively (“u” and “d” for short) after the two dominating computer operating systems which differ in their bent for the one or the other.

The name of the entity we’re focused on herein is suggested by the “42 Assessors” of the Egyptian Book of the ***Dead***, who sit in two rows of 21 along opposite walls of the Hall of Judgement.[5] During the soul’s so-called “Negative Confession” (and our focus is almost solely on units which are roots of negative unity), they ask embarrassing questions while the heart of the deceased is weighed in a ***pan-balance*** and the soul-eating Devourer of the ***Unjustified*** watches. The three senses of Zero implicit in the italicized words (and the need for two pans to achieve a balanced measurement), added to the coincidence of the count, made the name seem an appropriate mnemonic.

Our strategy, in some regards, is opposite Moreno’s. He approaches things head on by attempting to typify the non-standard or “special” features found in Sedenions that make zero-division possible. Here, we’ll be using everything we can that’s nascent in associative triplets – be they pure Octonions (“***O-trips***” henceforth) or inclusive of a pair of pure Sedenions (“***S-trips***”) – as a sprinter would exploit a thrust-block. Ditto for Moreno’s results. Rather than *concluding* with a “$G_2$” existence proof, we *begin* with a dissatis-fied need for a concrete revelation which the “$G_2$” finding has exacerbated.

Recalling that zeros of equations are roots, and $G_2$ ***root*** structure is an overlay of two hexagonally symmetric kaleidoscopes (a “doubled $A_2$”), this means two distinct orbits of differing “root length.” This latter designation makes sense in the standard 3-D renderings of the 12 roots of the algebra found, for instance, in Gilmore[6]. In the 8-D of an Octonion copy, the 6-cycle weaving between DOS and UNIX paths of Assessor triples, revealed by our first production rule, will provide the correspondence to the classic $A_2$ image. Our second production rule will show how any Assessor pair can spawn another, in a distinct $A_2$-like orbit, by a skew-symmetric swap of their non-O-trip units. Thanks to the Cartan subalge-bra, though, Lie algebraists think of $G_2$ as 12 + 2 = ***14***-dimensional, displayed most blatantly in the fact of the 7 O-trips, each equipped with a binary orientation (they’re associative, but ***not*** commutative!).

To keep the tools clean, we'll refer not to "$G_2$" but ***GoTo*** listings of the 12 Co-Assessor matches associated with each O-trip. Our third production rule, finally, will take us beyond *local* O-trip specifics to links *between* Octonion copies in Sedenion space (each Assessor appears in exactly ***2*** of the 7 GoTo lists).

We'll be short-handing zero-divisor pairs by use of their imaginary units' index numbers. Equations which would solve for these will be written as either (A, B) and (C, D) when signing doesn't matter, or else as either (A + B) or (C – D), to indicate a pair of partially or completely unknown units whose path of interest is UNIX or DOS respectively. (Hence, (– A – B) is UNIX, so (A + B) covers it, while (– C + D) is covered by ( C – D ).) Often, we'll have recourse to the binary variable "sg" (which can only be $\pm 1$) operating, say, between a '+' or '–' and the right-hand letter within the parentheses.

Most critically, we will indicate the product of two units A and B of unknown index like this: A *xor* B. This leads to the all-important convention in choice of labeling the Sedenions: the index of the product of two units must be the result of the XOR of the indices of the two units being multiplied. This requires some (brief) commentary.

The number of different ways one can label the Octonion units with suitable indices (480) is bad enough; but Sedenions are even worse (half a trillion!). Different authors have different tastes, and mine require products' indices be the XOR of the indices producing them. This is well-known to be possible, but seldom exploited as a key feature (but see Michael Gibbs' pages giving "Bit Representations of Clifford Algebras" posted on Tony Smith's web-site[7]). Without this simplifying constraint, in fact, none of the many "bit-twiddling" proofs to follow would be possible! Hence, units indexed 7 and 12 should have product with index (0111 XOR 1100 = 1011 =) 11.

It's impossible to make signing completely consistent with cyclic ordering, so I'd rather solve for than worry about it. I'm using the "canonical" Sedenion multiplication table, generated by unvarnished Cayley-Dickson process, found on Tony Smith's web-site (derived, in turn, from Lohmus, Paal and Sorgsepp's book[8]). As will become clear toward the end, this is no mere convenience. Where zero divisors are concerned, the assumed "equivalence" of the myriad multiplication tables allowed by the imaginary units' permutations and signings (already not quite perfect with the Octonions, as recent work by Schray and Manogue[9] has shown) breaks down. Whether this is "bug" or "feature," why it happens, and what to do about it, are all questions that will concern us later, but for which final answers are not yet at hand.

In the table, the real units are rendered as signed letter U's for purposes of this tabulating chore only. This much said, here's my Sedenion table, which readers with excess time on their hands should feel free to test for its XOR-worthiness:

| | **U** | **1** | **2** | **3** | **4** | **5** | **6** | **7** | **8** | **9** | **10** | **11** | **12** | **13** | **14** | **15** |
|---|---|---|---|---|---|---|---|---|---|---|---|---|---|---|---|---|
| **U** | U | 1 | 2 | 3 | 4 | 5 | 6 | 7 | 8 | 9 | 10 | 11 | 12 | 1 | 14 | 15 |
| **1** | 1 | -U | 3 | -2 | 5 | -4 | -7 | 6 | 9 | -8 | -11 | 10 | -13 | 12 | 15 | -14 |
| **2** | 2 | -3 | -U | 1 | 6 | 7 | -4 | -5 | 10 | 11 | -8 | -9 | -14 | -15 | 12 | 13 |
| **3** | 3 | 2 | -1 | -U | 7 | -6 | 5 | -4 | 11 | -10 | 9 | -8 | -15 | 14 | -13 | 12 |
| **4** | 4 | -5 | -6 | -7 | -U | 1 | 2 | 3 | 12 | 13 | 14 | 15 | -8 | -9 | -10 | -11 |
| **5** | 5 | 4 | -7 | 6 | -1 | -U | -3 | 2 | 13 | -12 | 15 | -14 | 9 | -8 | 11 | -10 |
| **6** | 6 | 7 | 4 | -5 | -2 | 3 | -U | -1 | 14 | -15 | -12 | 13 | 10 | -11 | -8 | 9 |
| **7** | 7 | -6 | 5 | 4 | -3 | -2 | 1 | -U | 15 | 14 | -13 | -12 | 11 | 10 | -9 | -8 |
| **8** | 8 | -9 | -10 | -11 | -12 | -13 | -14 | -15 | -U | 1 | 2 | 3 | 4 | 5 | 6 | 7 |
| **9** | 9 | 8 | -11 | 10 | -13 | 12 | 15 | -14 | -1 | -U | -3 | 2 | -5 | 4 | 7 | -6 |
| **10** | 10 | 11 | 8 | -9 | -14 | -15 | 12 | 13 | -2 | 3 | -U | -1 | -6 | -7 | 4 | 5 |
| **11** | 11 | -10 | 9 | 8 | -15 | 14 | -13 | 12 | -3 | -2 | 1 | -U | -7 | 6 | -5 | 4 |
| **12** | 12 | 13 | 14 | 15 | 8 | -9 | -10 | -11 | -4 | 5 | 6 | 7 | -U | -1 | -2 | -3 |
| **13** | 13 | -12 | 15 | -14 | 9 | 8 | 11 | -10 | -5 | -4 | 7 | -6 | 1 | -U | 3 | -2 |
| **14** | 14 | -15 | -12 | 13 | 10 | -11 | 8 | 9 | -6 | -7 | -4 | 5 | 2 | -3 | -U | 1 |
| **15** | 15 | 14 | -13 | -12 | 11 | 10 | -9 | 8 | -7 | 6 | -5 | -4 | 3 | 2 | -1 | -U |

For $2^N$-dimensional "blow-ups" of the Imaginaries up through the arithmetic of interest, each unit, once the Reals and the unit itself are excluded from the count, is involved in $(2^N - 1)(2^N - 2)/(3!)$ associative triples. There's but one such triple for the Quaternions, 7 for the Octonions, and 35 for the Sedenions. These are the staff of life for all that follows, so let's stop to get acquainted with them.

Let's short-hand associative triples as follows: if (ab)c = a(bc), take the index numbers of these imaginary units and write them as a comma-delimited list in parentheses, in cyclical order. The lowest index number comes first, then the index whose product with it yields the positive value of the third. The Quaternions' "i, j, k" – almost always written as $i_1$, $i_2$, and $i_3$ respectively – has its associative triplicity indicated this way: (1, 2, 3). Using our Sedenion table, we can write the seven O-trips as follows:

(1, 2, 3); (1, 4, 5); (1, 7, 6); (2, 4, 6); (2, 5, 7); (3, 4, 7); (3, 6, 5)

It's visually obvious that only two of these triplets have counting order varying from cyclical order: $i_1 \times i_6 = (-i_7)$, and $i_3 \times i_5 = (-i_6)$. (Of course, different triples are out of order in different notations.) The other 28 associative triples in the Sedenions – the S-trips – are listed, per our table, next:

| | | | |
|---|---|---|---|
| (1, 8, 9) | (1, 11, 10) | (1, 13, 12) | (1, 14, 15) |
| (2, 8, 10) | (2, 9, 11) | (2, 14, 12) | (2, 15, 13) |
| (3, 8, 11) | (3, 10, 9) | (3, 15, 12) | (3, 13, 14) |
| (4, 8, 12) | (4, 9, 13) | (4, 10, 14) | (4, 11, 15) |
| (5, 8, 13) | (5, 12, 9) | (5, 10, 15) | (5, 14, 11) |
| (6, 8, 14) | (6, 15, 9) | (6, 12, 10) | (6, 11, 13) |
| (7, 8, 15) | (7, 9, 14) | (7, 13, 10) | (7, 12, 11) |

The leftmost diagram at the bottom of the next page is the "canonical" labeling of PSL(2,7) resulting from applying the Cayley-Dickson process. Each projective line contains 3 points (which determine the "circle" of the Octonion 3-cycle indicated by the points' labels). The cycle in the middle (2, 4, 6) is the only one actually appearing in this representation *as* a circle, although it is a "line" just like the other six. This much is true for *all* representations employing the triangle. The differences concern labels and arrows only: some authors also use an XOR representation, and so keep the labeling, but change the sign of a unit, and so change the directions of the three arrows that pass through it. Or, they change multiple units' signs, thereby switching some other count of arrows. (But always all or none, or 3 or 4, as we'll see later on).

The most frequently encountered variant (used by Schray and Manogue) switches the 6's sign, and hence reverses the ( 1 7 6 ), ( 3 6 5), ( 2 4 6 ) triplets' cycling direction, replacing them, per our earlier convention, with ( 1 6 7), ( 3 5 6 ), and ( 2 6 4 ) respectively. This shift can be represented by simply flipping the bottom, vertical, and circular arrows; or (as Schray and Manogue do), by rotating the two lines that pass from 1 through 6 and 5, which thereby also reverses the arrows on the diagonals (but now labeled 3, 4, 7 and 6, 4, 2 in ascending order) as a side-effect. Then, change the orientation on the circle in the middle and the right-hand slope to have labels and arrows properly arranged.

The appeal of this variant is that it offers the minimum deviation from counting order in cyclically ordered triplets: only one ( 2 6 4 ) is out of order in this labeling. But more extreme variations can also be found, which are strictly combinatoric in nature. Feza Gursey, a pioneer in the study of "theory of every-thing" models employing the exceptional Lie algebras like $E_6$, used a labeling system which didn't display XORing, but did offer some analogy to tensor-indexing schemes general relativity theorists might use. Other physicists focused on model-building with unusual algebras have also used this scheme, or – like Susumu Okubo[10] – a sign-flipped variant (the 7 in Gursey's version has its sign reversed, yielding the diagram in the middle at the bottom of this page.)

Some recent work of great interest has focused on exploiting the 2- and 3-term cross-products which are unique to the Octonions. (The familiar cross-product based on the Quaternions, and adapted since Gibbs to "vector mechanics," has no other true analogues[11], and they also have a "$G_2$" symmetry.)

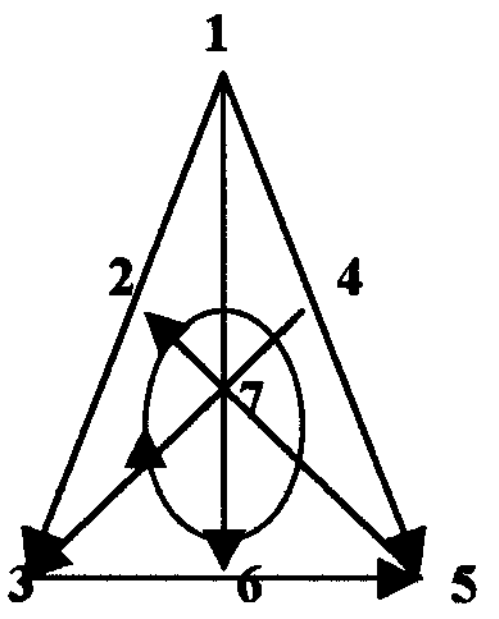

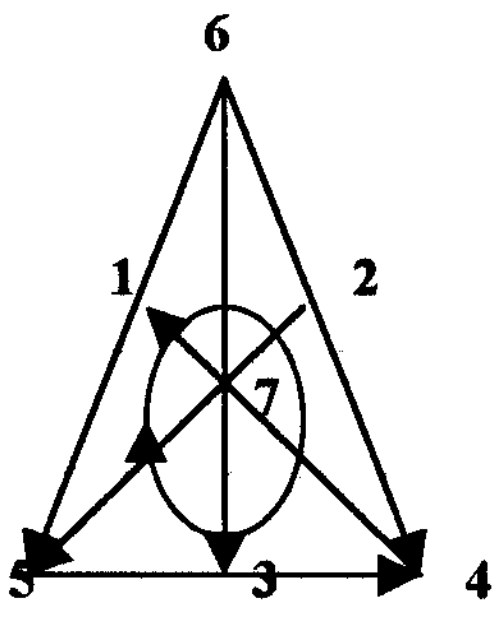

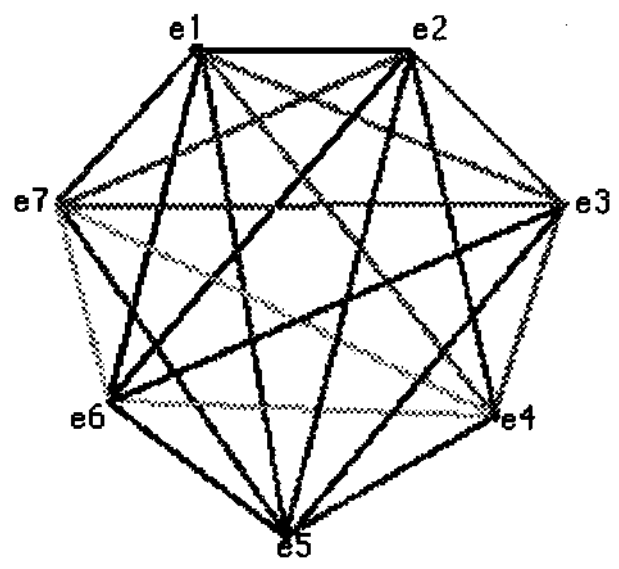

Onar Aam and Geoffrey Dixon have been especially fond of numbering schemes based on mod-ulo-7 arithmetic, rather than XORing, and hence don't even use the PSL(2,7) representation. Instead, they draw triangular circuits (one per O-trip) connecting the vertices of a heptagon, as shown above.

The cyclicality of indices seems the major draw here. With modulo 7 arithmetic, one can avail oneself of more than one equivalent way of representing this: $e_i * e_{i+1} = e_{i+5}$ is isomorphic to $e_i * e_{i+2} = e_{i+3}$, and mirror-opposite (through the vertical) to $e_i * e_{i+1} = e_3$ ( or $e_i * e_{i+2} = e_{i+6}$ ). By developing special "X," "XY," and "WXY" products based upon the Octonions' cross-product structures, Dixon[12] has used his heptagonal indices in deep studies of 16- and 24-D exotica of great interest to superstring theorists, the closest-packing-pattern lattices of Barnes-Wall and Leech respectively. Aam, meanwhile, has done fascinating work on Octonion cross-product symmetries, leading him to construct things Tony Smith calls "onarhedra" – and which, I discovered near the end of writing this monograph, are 8-D isomorphs of the "box-kites" I'd independently derived from 16-D configurations of zero-divisors. (After reading this paper, a visit to Smith's web-site is highly recommended, where the pages on Aam's "cross-product" work will show amazing correspondences – once one adjusts to notational and intuitional differences – to the entities I've called "strut signatures" and "tray-racks," for instance. All of which gives further confirmation of the profundity, qua derivation algebra at least, of the "$G_2$ connection" between Octonions and Sedenions.)

The 480 different multiplication tables usually thought to be equivalent ways of working with Octonions may seem like a lot – until one computes how many one must deal with in higher-dimensional cases. For the Octonions, the 30 distinct combinatorial variations (up to signing) result from dividing 7! = 168 x 30 by the internal symmetry group of the finite projective plane whose labels and arrows we've played with – the 168-element simple group, PSL(2,7). With the Sedenions, we have to divide 15! ~ $1.24 \times 10^{15}$ by PGL(4,2) = 8!/2 = 20,160 (the symmetry group of finite projective *space*), then multiply *that* by the number of sign-flip-only possibilities: 2,048 = $2^{15 - 4}$ (4 = number of generators via Cayley-Dickson process, 3 for Octonions). The result is half a trillion . . . and these, as we'll see, are *not* equivalent!

The Sedenions' indices can be placed on the 15 vertices of the finite projective *tetrahedron* with 35 lines of 3 points each, one per triplet, 7 lines intersecting in each point and 3 planes in each line, implying 1 Octonion triangle + its 7 lines x (2 – 1) distinct planes through each = 15 planes isomorphic to the PSL(2,7) triangle. The picture is from p. 139 of Lohmus et al., where it found its way to Tony Smith's web-site, where *I* first encountered it. But mostly, working with triplets directly will seem far more natural and easy.

| | |
|---|---|
| **Real numbers**<br>– the corresponding projective geometry consists of empty set of points (dim =−1) | |
| **Complex numbers**<br>– the corrresponding projective geometry consists of a single point (dim = 0) | ∘ $e_1$ |
| **Quaternions**<br>– the corresponding projective geometry consists of three collinear points (dim = 1) | 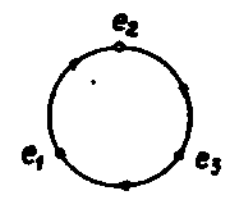 |
| **Octonions**<br>– the corresponding projective configuration consists of 7 points and 7 lines (3 points on each) (dim = 2) | 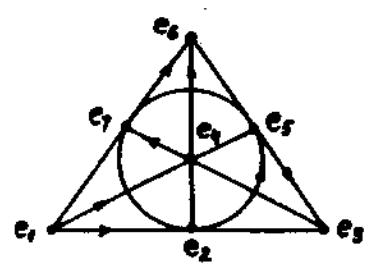 |
| **Binary sedenions**<br>– the corresponding projective 3-dimensional configuration consists of 15 points and 35 lines, 7 lines through every point; from purely technical reasons 18 lines are missing on figure | 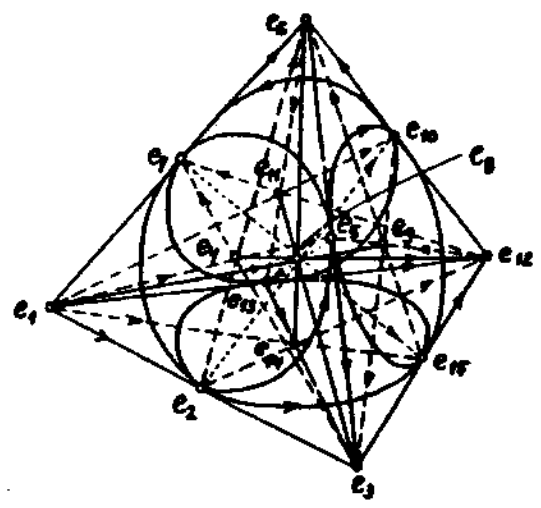 |

While Moreno generates his $G_2$ analysis in the usual "top-down" manner, giving but one instance of an actual zero-divisor pair in his paper, my own bias is toward "bottom-up" constructivism. The process of creating actual zero-divisor pairs will occupy us next. Our short-term goal will be to provide all things necessary for revealing the 42 Assessors, proving their completeness, and generating "GoTo lists" for the seven O-trips which will provide necessary and sufficient kindling to spark all this. (We'll see, too, that all but seven of the S-trips are included in passing – three in fact, per O-trip GoTo list. The exclusion of seven others from all such lists will be proven by our "Behind the Eight-Ball Theorem" late in our construction work, and will motivate another tabulation, the "***Osiris Partition***," from which ***box kites*** will fly.)

## { II } Creating Something from Nothing (Three Production Rules):

Quartets of imaginary units written per the notational conventions given above form valid zero-dividing pairs if the following obtains:

(A, B) x (C, D) = 0 ⇔ (A *xor* C) – (D *xor* B) = (A *xor* D) – (C *xor* B) = 0.

Alternatively, given the correlation between signing and left-right orientation as one moves in sequence in an associative triplet's 3-cycle, we could just as readily write the above rule this way:

(A, B) x (C, D) = 0 ⇔ (A *xor* C) + (B *xor* D) = (A *xor* D) + (B *xor* C) = 0.

In the realm of the Octonions, only trivial cases entailing null or non-distinct letters can occur. This is certainly connected to the fact that each of the seven Octonion triplets has exactly one unit in common with any other triplet. Translating this into conditions among Sedenions is not immediately obvious, however. First, any product like the above implies *six* different triplets, since each Sedenion index pair belongs to the triplet containing its XOR product. Hence, A and B form a triplet with A *xor* B = X, while C and D form another with C *xor* D = Y. Moreover, the requirements just given mandate that *two pairs* among the six must necessarily share a term each: (A, C) and (B, D) share some E; (A, D) and (B, C) share some F.

In fact, we can say more than this. Since the Octonions and their copies contain but seven imaginary units, and the explicit requirements for zero-dividing already exhaust six, then:

<u>Lemma:</u> Given two Co-Assessors (A, B) and (C, D), A *xor* B and C *xor* D have the same index G.

This clearly implies that the units involved in multiplication can be arranged in both Assessors so that (A *xor* B ) - ( C *xor* D) = 0, without effect on the end result. This "axiom of choice" will be of great use in the proof behind our second production rule. For the moment, though, we provide, prove, and apply

<u>Production Rule #1 ("Three-Ring Circuits")</u>. Given a pair of Co-Assessors (A, B) and (C, D), a new Assessor ( E, F ), mutually zero-dividing with both its progenitors, can be created by assigning the new letters to the two distinct XOR's ( A *xor* D = C *xor* B; A *xor* C = D *xor* B ).

<u>Proof:</u> Do the zero-dividing arithmetic, keeping track of the triplet cycles as they reveal themselves.

C + D

A + B

[B *xor* C = ] + F  + E [ = B *xor* D ]  => ( B, C, F ); ( B, D, E )

[A *xor* C = ] – E  – F[ = A *xor* D ]  => ( D, A, F ); ( C, A, E )

0

Now, substituting from the associative triplets listed to the right, form products of ( E ± F ) with the initiating Co-Assessors, allowing for unknown signing via the “sg” variable:

| | |
|---|---|
| E + sg F | E + sg F |
| A + B | C + D |
| - D - sg C | + B - sg A |
| + C + sg D | - A + sg B |
| 0 ⇔ sg = (+1) | 0 ⇔ sg = (-1) |

Clearly, exactly one of these must be true for fixed sg. By moving sg between lines of a product, it is obvious that we get a collection of six distinct product-pairs involving ( A ± B ), ( C ± D) , and ( E ± F ). As we make zero by multiplying the Assessors in some sequence, we are forced to switch from DOS to UNIX paths in one of two cross-over sequences, the knot theorist’s trefoil (dduddu), or else the “unknot” of a “triple zigzag” (dududu). In the GoTo Listings, the first of each O-trip’s 4 columns displays, in 3 consecutive lines, the first half of its sole zigzag’s 6-cycle; the other 3 columns show the 3 trefoil couplings.

Dynamic Content. Whether trefoil or zigzag, the dynamics are analogous – and require adapting the usual exponential power-orbit to the special requirements of their nonassociative (and even non-alternative) world. The complex domain’s exp(ix) = cos x + i*sin y requires various “correction terms” to work in less orderly settings. An Assessor’s pair of imaginaries indexed ( A ± B ) can’t simply be recast as exp( $i_A$ ) (cos x + $i_C$ sin y), A *xor* C = ± B, with a Co-Assessor indexed ( D ± E ) being recast similarly. For the zero to show in the product, the sine and cosine terms in the right-hand parentheses must evaluate ***last***.

The simplest way to guarantee this is by fiat: define a "zip" (Zero-divisor Indigenous Power-orbit) function. Ignoring internal signing for simplicity, for Assessor indices ( A, B ) and ( C, D ), and angular variables x and y, write products of the circular motions in the two Assessor planes like this:

$$\text{zip}( A, B;\ x ) * \text{zip}( C, D; y ) = ( \cos x * i_A + \sin x * i_B ) * ( \cos y * i_C + \sin y * i_D )$$

$$= ( \cos ( x + y ) * i_{[A\ xor\ C]} ) + ( \sin ( x - y ) * i_{[D\ xor\ A]} )$$

We know, from our first production rule, how to interpret this: two currents with opposite clockwise senses or "chirality," spawned by angular "currents" in two Co-Assessors' planes, manifest as alternately con- and de- structive wave interference effects, in the plane spanned by the third Assessor they mutually zero-divide. In this "toy model," patterns readily recognizable as Lissajous figures would sweep through the origin, "showing on the oscilloscope" in accordance with the opposing currents' relative velocities (provided, of course, their ratio only involved small integers).

This is, of course, vastly simplified: just one Co-Assessor trio's interactions require considering six dimensions of toroidal "plumbing," and Assessors can interact directly in two different trios, and indirectly receive or transmit "pass-along currents" in ways that will require separate study. By the end of the present investigation, we will have sufficient information to enable the framing of such studies, but insufficient space to contemplate them further in these pages.

One point, though, must be made now: thanks to Peixoto's Theorem[13] (well known to dynamicists, not well known to high-energy physicists), we know that stable flow-patterns on a torus can only be generated by pairs of oppositely oriented currents with same rational winding number – clearly suggesting that zero-divisor dynamics in the badlands of Sedenion space could provide the basis for a surprising new kind of orderliness! (Consider the unpromising mix of dispersion and nonlinearity, yielding the supreme stability of soliton waves.) To see how naturally this condition can be met, let's make our "toy model" a bit more robust. Show the full 6-cycle, assume noncommutativity implies two helically distinct routes through the cycle, and think of them as generated by a "volley and serve" logic. (That is, the helices "converse," with the (A,B)-plane-initiated "ping" inducing a (C,D)-plane-based "pong" in some "dt" less than the time it takes for one full wind around the torus.) Here's "dt"-free for-instances of what this could look like:

| Ping | Pong |
|---|---|
| $+\cos(x + y)\, i_1 + \sin(x - y)\, i_{13}$ | $-\cos(x + y)\, i_1 - \sin(x - y)\, i_{13}$ |
| $+\cos(y + z)\, i_2 - \sin(y - z)\, i_{14}$ | $-\cos(y + z)\, i_2 + \sin(y - z)\, i_{14}$ |
| $+\cos(z + x)\, i_3 + \sin(z - x)\, i_{15}$ | $-\cos(z + x)\, i_3 - \sin(z - x)\, i_{15}$ |
| $+\cos(x + y)\, i_1 - \sin(x - y)\, i_{13}$ | $-\cos(x + y)\, i_1 + \sin(x - y)\, i_{13}$ |
| $+\cos(y + z)\, i_2 + \sin(y - z)\, i_{14}$ | $-\cos(y + z)\, i_2 - \sin(y - z)\, i_{14}$ |
| $+\cos(z + x)\, i_3 - \sin(z - x)\, i_{15}$ | $-\cos(z + x)\, i_3 + \sin(z - x)\, i_{15}$ |

It must be noted, as a caveat, that generalizations of Peixoto's Theorem beyond the 2-torus are nontrivial, and even problematic in many contexts. But the possibilities in the present case seem so beautiful that, barring solid evidence to the contrary, pessimism feels too difficult to sustain.

Two final notes: the ping/pong "dt" is the natural condition for allowing quantization, and as long as we expect responses to multiplicative acts are not immediate, we can readily construct toroidal flows in a dynamic "beneath" quantization which are stable in at least the limited 2-torus context. Later, we'll see that there are 28 "CoAssessor Trios," which suggests a link to the '28' of the D4 Lie algebra which, via "triality," can underwrite an enormous array of quantized structures. In the "hyperdiamond lattice" approach – an alternative or complement to the "superstring" theoretics that garner most of the press, linked with the work of David Finkelstein and Tony Smith and their students – the D4 focus is paramount in many respects. And, as one can readily discover by browsing his website, Smith has suggested its connection to the 28 "S-trips" in the Sedenions. I would suggest here that a deeper connection to the Sedenions may be found through the CoAssessor Trios – which brings up my second point.

The "triality" among zero-divisors forming a Trio has crucial elements not like the 2-torus at all – or at least not like it in any obvious way. What family of functions A, B, C and/or D, E, F can effect some kind of circuit in a Trio's 6-cycle, returning to some analog of "start-up mode"? Or (analogous to Fourier series), can *create* a cycle which will behave *as if* we'd had it as our starting point, irrespective of what sort of messiness the resonance emerged from? What *fourfold* function families can be envisioned for circuits involving *quartets* of zero-divisors (which we'll be calling "tray-racks" later)? Can some dense set of such constructs be built upon a Peixoto-style "bedrock" of 2-torus stable flows, or are things far wackier than this? Do we need a dynamic model based on "quick-hit" interventions and stably regulated jumps, akin to the needle-sticking redirecting of the flow of *chi* in acupuncture?

*    *    *    *    *

Production Rule #2 ("Skew-Symmetric Twisting"). Given Co-Assessors (A, B) and (C, D), precisely two new Assessors, unable to be Co-Assessors with either propagator but Co-Assessors with each other, can be created by exchanging one lettered unit from within one Assessor's parentheses with one from the other, then flipping the sign associated with one or the other unit in the resultant.

Preliminary Remarks. Since exchanging A and D is effectively indistinguishable from exchanging B and C, the four possible pair-swaps reduce to but two. Without loss of generality, we can then say that the restriction to "precisely one" means either ( A + D )( C – B ) or ( D + B )( C – A ) = 0, but not both.

Proof. By the earlier lemma, we can arrange the writing order of the terms in the parentheses of the two Co-Assessors we apply the rule to so that the XOR's of the terms internal to the parentheses have opposite sign. That is, for ( A + B )( C + D ), we arrange our labeling so that A *xor* B = + G, and C *xor* D = - G. This gives O*-trip sequences ( A, B, G ) and ( D, C, G ) respectively. (Note: write O*- and not O-trip, as we are necessarily dealing with an Octonion *copy*! All but *one* of the O*-trips within it will, of necessity, be S-trips and hence *not* be contained within the list of seven O-trips proper given earlier. More on this later!)

Maintaining the same notation otherwise as with our prior proof for Production Rule #1, yields the other O*-trip sequences: ( B, C, F ) and ( D, A, F ); and ( B, D, E ) and ( C, A, E ). Let's do the arithmetic:

| | |
|---|---|
| C - B | C - A |
| A + D | D + B |
| + G + E | + F + G |
| - E - G | + G - F |
| 0 | 2G |

Clearly, the left-hand product is generated by a new pair of Assessors. Just as clearly, products of either with either generator will result in real units, since each got one letter each from each propagator. In normed units, the "2G" result in the right-hand product reduces in fact to plus or minus the normed imaginary unit produced by the internal XOR-ing described above.

These two production rules can be combined in two ways. First, Rule #1 can be run on the results of Rule #2. But also, Rule #2 can be run on the results of Rule #1! This means successive applications of the two rules turn two Assessors into twelve . . . a result strong enough to require a

Proof. No checking is needed for each Assessor pair before the Twisting and that created after it. We need, though, to check that results of the three Twistings do not duplicate each other ( 3 x 2 = 6 tests), and do not recreate the third Assessor related to each Twisting's propagating pair by Production Rule #1. By shared letters, the six Assessors created by Twisting, can each only require testing for distinctness against one of each of the two other Twistings' Co-Assessor pairs, which means ( 3 x 4 )/2 = 6 separate tests. So 12 tests in all are called for to complete the proof, which symmetry can reduce to 6. Consider the Twisting of ( A + B ) with ( E – F), which produces ( A - F ) and ( E - B ). Test both against ( C + D ); ( A - F ) against ( C – B ) and ( E + D ); and, ( E - B ) against ( A + D ) and ( C + F ).

We choose this Twisting to examine, because it completes the full mapping of the Octonion copy by resolving the seventh O*-trip we know must exist in such a copy. We know ( A, B, G ), and this Assessor pair requires this to be negated by E *xor* F. Therefore, the seventh O*-trip reads: ( E, G, F ). For convenience, all seven O*-trips in our O-copy are collected together now:

( A, B, G ); ( A, E, C ); ( A, F, D ); ( B, C, F ); ( B, D, E ); ( C, G, D ); ( E, G, F )

Here's a sequence of six quick calculations to conclude the proof:

( A – F )( C + D ) = ( – E – F + B – A ) FAILS! ( E – B )( C + D ) = ( + A – B – F – E ) FAILS!

( A – F )( C – B ) = ( – E – G + B + C ) FAILS! ( A – F )( E + D ) = ( + C – F – G – A ) FAILS!

( E – B )( A+ D ) = ( – C – B + G – E ) FAILS! ( E – B )( C + F ) = ( + A – G – F + C ) FAILS!

Masochists who won't take "arguments by symmetry" at face value are encouraged to spend a minute or two crunching out three more lines of calculations. They may also enjoy playing with the four-unit products resulting from all failed results. Each includes precisely two terms from its producers, and, with norming assumed, they are all in fact square roots of minus unity, as can be checked quite quickly. (All cross-products among the four letters cancel, the squares of the units sum to –4 within the parentheses, and the norming of the producing terms mandates the coefficient of the parenthesized listings above be ½.)

After a little playing around like this, though, an earlier failure should return to your thoughts, which also resulted in an Imaginary (instead of pregnantly Empty) result. What was it about the *second* possible Twist result that made it fail? This is related, of course, to our exercising the free choice our earlier lemma showed was our prerogative. And perhaps it is related to another fact tossed off in a parenthetical remark a few lines back. Of the seven O*-trips in an Octonion copy, only (and *exactly)* one can actually be an O-trip proper . . . belong, that is, to the genuine Octonions (index $< 8$) prior to their extension to Sedenions.

Let's mull on this for a minute, beginning our thinking with a new convention. To always be clear which units have index > 7 (pure Sedenions), we will use letter variables that are Upper Case only, such as 'S'; correspondingly, pure Octonions will be written exclusively with lower case lettering ('o'). No copy O* of the genuine Octonions can be all S's, for the XOR of any two will produce a Zero in the leftmost bit, and hence result in an o. Ditto, all o's means the genuine O, hence no division by Zero. But one of each mandates *two* of each by XORing (o *xor* S = some S'; S *xor* S' = some o'), and two o's have a third as product, completing an O-trip.

The third o, when XORing with S and S', must produce at least one new S, call it S''. Suppose the first o forms an associative triple with S and S'. Then either o' or o'' must form an S-trip with two S's distinct from both of these. This leaves us with 3 o's and 4 S's, which completes the O*. But suppose o' does *not* form an associative triple with S and S'. Then either o' or o'' must form an S-trip with them instead, which leads to the same result. Hence, any O* not identically the genuine O must have 3 o's and 4 S's exactly.

Our Co-Assessor pairs explicitly require writing with four letters. All four cannot be o's, lest we be in the genuine O. We cannot have one o or S, and three S's or o's, lest multiplying require an S *xor* o to equal an S *xor* S or an o *xor* o, which is impossible if we're to end up with Zero. So either *all* must be S's (which we'll see is impossible); or, we must have two of each.

Suppose the latter is the case. Then we can have it two ways: (1) the o's are combined to make one Assessor, and the S's, the other; or, (2) an Assessor must have one of each. Regardless of which case we *start* with, the two possible Twists the proof of Production Rule #2 required we test (swapping terms on the same or opposite sides of the middle + or -) will *result* in exactly one case of each Co-Assessor pairing. And we know that *exactly one* of their multiplications must *always* fail. Which suggests the

“Exogamous Moieties” Rule. It is forbidden for the two letters on either side of an Assessor’s central ‘+’ or ‘-’ to belong to the only O-trip among an O-copy’s seven associative triplets. Put another way, the two units whose marriage constitutes an Assessor must belong to different halves or “moieties”of the Sedenion tribe. One must be a pure Octonion or ‘o’; the other must be a pure Sedenion or ‘S’.

Proof. Consider the ways of picking two units. If moiety didn’t matter, then other things being equal, there should be (7x6)/(2x1) = 21 possibilities. Picking 2 of 3 o’s can happen 3 ways; picking 2 of 4 S’s, (4x3)/(1x2) = 6 ways; and picking one of each can be accomplished 3x4 = 12 ways, which three options tally, of course, to 21. But we have found precisely 12 Assessors linked to any given O-trip. Moreover, without even knowing actual index numbers for *any* actual cases yet, we can write them out abstractly and catch the pattern anyway. Writing each Co-Assessor triple in a separate column, we have this array:

| | | | |
|---|---|---|---|
| ( A, B ) | ( A, G ) | ( A, D ) | ( A, F ) |
| ( C, D ) | ( C, F ) | ( C, B ) | ( C, G ) |
| ( E, F ) | ( E, D ) | ( E, G ) | ( E, B ) |

If our assertion be true, then A, C, E are an O-trip (and should be written “a, c, e” by our casing convention). In the first test of Production Rule #2, it was ( A, C ) ( D, B ) which failed. Rewritten to show casing, this is ( a, c ) ( D, B ), confirming our hypothesis. Moreover, twisting ( A, B ) ( E, F ) fails for the ( a, e ) ( B, F ) product only; likewise, ( C, D ) ( E, F ) has ( c, e ) ( D, F ) as a failed twist product. And easy calculations show this obtains with the three other columns as well.

The last possibility for failure of the proof is the four-S case, which we’d placed to the side. And however we form an instance of it, *all four XOR products must be o’s*. But then “Exogamous Moieties” must be invoked, and so the First Production Rule must fail, which contradicts the claim we have a Co-Assessor pair in the first place. Q.E.D.

Afterthoughts: How To Instantiate and Interpret an “O-Copy”: The “letter logic” of the just concluded argument is not quite a concrete tool for producing a verifiable “O-Copy.” For that, we need to massage the labels by making use of the last theorem in Moreno's paper (the solitary corollary of which is his "$G_2$ homomorphism" result). Put simply, he proves that if we have two units a, b which associate, and a third we label y which anti-associates with them – ( ay )b = – a( yb ), in other words – then a, b and y provide necessary and sufficient material for spawning a copy of the Octonions. A mapping from (a,y,b)-generated terms to the Octonion index set is made by this sequence of assignments:

( a, b, ab, (ay)b, yb, ay, y ) ➔ ( 1, 2, 3, 4, 5, 6, 7 )

Let's set the O-trip ( A E C ) from our proof apparatus to be ( a, b, ab ); and, let's pick D for our y. Plugging in the rest of our letters into the mapping machinery gives us these assignments:

( a, b, ab, (ay)b, yb, ay, y ) ➔ ( A, E, C, -G, B, -F, D ) ➔ ( 1, 2, 3, 4, 5, 6, 7 )

To give this a test drive, let's stuff some numbers. Pick ( a, b, c ) to be the Quaternions' ( 1, 2, 3 ) and then search the GoTo Listings until we get a Co-Assessor Trio associated with this O-trip which induces minimal signing messiness when we run it through the steps of our proofs. Here's the result:

( 1 + 14 )( 3 + 12 ) = ( 3 + 12 )( 2 – 13 ) = 0 ➔ ( A + B )( C + D ) = ( C + D )( E – F ) = 0

All that's missing is G, but since AB = FE = G, G must be 15. That means we get this map:

( 1, 2, 3, –15, 14, –13, 12 ) ➔ ( 1, 2, 3, 4, 5, 6, 7 )

Since all GoTo Listings for a given O-trip share the same quartet of Sedenion units, the "O-copies" associated with each of the four O-trip's four Co-Assessor trios will employ the same (albeit differently signed) units in their mappings. Hence the Listings refer not to an automorphism, but an "automor***pheme***" (e.g., set of "constituent units," as they say in structural linguistics, from which an actual Co-Assessor trio's automorphisms are generated) when listing the seven imaginaries from which O-copies are made.

There is a secondary reason for introducing this subtle distinction: while we do, in fact, have copies of the Octonions, and hence O-copies in the sense used above, the "isomorphism" Moreno depends upon in fact entails a crucial loss of information. If we take the mapping just provided and generate the constituent O-trips for the terms on the right of the arrow, we'll see they don't quite match the S-trips which spawned the listing on the left. We have, in fact, ***two signing conflicts***.

( 3, 4, 7 ) maps to ( 3, –15, 12 ) ~ ( 3, 12, 15 ). Our Sedenion table shows the S-trip as ( 3, 15, 12 ).

( 3, 6, 5 ) maps to ( 3, –13, 14 ) ~ ( 3, 14, 13 ). Our Sedenion table shows the S-trip as ( 3, 13, 14 ).

This is a possibility which Moreno apparently didn't consider: indeed, the only concrete instance he offers of his "special triples" theorem is, perversely, one which has no zero-divisors: a, b, y = 1, 2, 7 . . . which generates not an O-copy, but the Octonions proper! Moreover, the signing conflict cannot be removed by permuting the signs on the underlying units. This is best checked on the Octonion side of the mapping, as there are only 480 equivalent multiplication tables. These partition into 30 combinatorially different placements on the PSL(2,7) triangle with same signing, times 16 identically labeled, differently signed, placements. The XOR requirement made at the onset restricts us completely as to combinatorial varying. We are, therefore, restricted to the 16 signing variants. But these in fact can be reduced by half, since each multiplication table has an "opposite" table with all arrows reversed.

We are down, then, to eight variants, each of which is generated by changing the sign on at least one unit, and hence the orientation of the 3 arrows which pass through it. But we needn't look at any of them explicitly, since the fact that any two O-trips share precisely one unit is all we need to know to prove this

<u>Lemma.</u> Permuting signs on units changes 3 or 4 arrows, or all or none.

<u>Proof.</u> Changing one unit's sign flips 3 arrows. Changing any two flips 6 – 2 = 4 arrows (since sign-change is an involution, the O-trip with both units is left unchanged). Changing all 3 members of an O-trip changes *all* arrows (the O-trip in question is flipped thrice = once; all other O-trips share but one, so get one flip). Changing all 4 units *not* in some O-trip changes *no* arrows (the excluded O-trip is left untouched, and all others are toggled twice). Changing 3 ( = 4 – 1 ) units which are not an O-trip therefore changes 3. By duality of opposite multiplication tables (there being but 7 units), this exhausts all possibilities.

<u>Throwaway Corollary.</u> Since the "canonical" Octonion multiplication table has two O-trips with cyclical and counting order differing, a table with all O-trips in order in both senses is impossible.

Implication. Moreno's "isomorphism" between Sedenion subsets which act like O-copies in their Sedenion context, and systems of units which are multiplicatively equivalent to the Octonions, is only partial, and breaks down where signing of units (and hence orienting of O-trips) of a true Octonion multiplication table is concerned. Moreno's "isomorphism" also doesn't capture sign-related information that's clearly germane to Co-Assessor Trios. ("Homomorphism to $G_2$," where roots or "zeros" are concerned, says nothing about the sign-based distinction between zigzags and trefoils: it merely points to the existence of "distinct" 6-cycle orbits of some kind – and it's even ambiguous about what "distinct" means, since not just our second production rule, but the third we'll see next, could claim, for different reasons, to define it.)

Moral: Tools derived from Lie algebra shouldn't be assumed to work as planned with things that are *not* Lie algebras. When determining laws about apples from the study of oranges, concrete data should at least be looked at before claiming they fit the facts! (Else there may be hell to pay later, as we'll find out!) This much said, here, then, are some further facts: the two mis-signed triplets are not due to a "mistake" by Moreno, they also are built in to the abstract "letter logic" of our proofs (readers are encouraged to test this out for themselves). More peculiar still, all the mis-signed triplets are the *same* triplets, regardless of what Sedenions or O-trips are used to feed Moreno's mapping machine. Further, all excluded Sedenions (and jump ahead to the "8-Ball Theorem" to see just what this means) generate "normal" Octonion multiplication tables. These are computable *facts*, whose underlying reasons we leave, for now, as mysteries.

* * * * *

## Table of the 42 Assessors by GoTo Listing, Each Appearing in Two, With Zero-Dividing Partners Shown Coupled, and Trios Shown in Columns

**GoTo #1 Based on Octonion Triplet (1, 2, 3) Automorpheme: (1, 2, 3, 12, 13, 14, 15)**

| | | | |
|---|---|---|---|
| (1 + 13)(2 – 14) | (1 + 14)(2 + 13) | (1 – 12)(2 – 15) | (1 – 15)(2 + 12) |
| (2 – 14)(3 + 15) | (2 + 13)(3 – 12) | (2 – 15)(3 – 14) | (2 + 12)(3 + 13) |
| (3 + 15)(1 – 13) | (3 – 12)(1 – 14) | (3 – 14)(1 + 12) | (3 + 13)(1 + 15) |

**GoTo #2 Based on Octonion Triplet (1, 4, 5) Automorpheme: (1, 4, 5, 10, 11, 14, 15)**

| | | | |
|---|---|---|---|
| (1 + 14)(4 – 11) | (1 + 11)(4 + 14) | (1 – 15)(4 – 10) | (1 – 10)(4 + 15) |
| (4 – 11)(5 + 10) | (4 + 14)(5 – 15) | (4 – 10)(5 – 11) | (4 + 15)(5 + 14) |
| (5 + 10)(1 – 14) | (5 – 15)(1 – 11) | (5 – 11)(1 + 15) | (5 + 14)(1 + 10) |

**GoTo #3 Based on Octonion Triplet (1, 7, 6) Automorpheme: (1, 7, 6, 10, 11, 12, 13)**

| | | | |
|---|---|---|---|
| (1 + 11)(7 – 13) | (1 + 13)(7 + 11) | (1 – 10)(7 – 12) | (1 – 12)(7 + 10) |
| (7 – 13)(6 + 12) | (7 + 11)(6 – 10) | (7 – 12)(6 – 13) | (7 + 10)(6 + 11) |
| (6 + 12)(1 – 11) | (6 – 10)(1 – 13) | (6 – 13)(1 + 10) | (6 + 11)(1 + 12) |

**GoTo #4 Based on Octonion Triplet (2, 4, 6) Automorpheme: (2, 4, 6, 9, 11, 13, 15)**

| | | | |
|---|---|---|---|
| (2 + 15)(4 – 9) | (2 + 9)( 4 + 15) | (2 – 13)(4 – 11) | (2 – 11)(4 + 13) |
| (4 – 9)(6 + 11) | (4 + 15)(6 – 13) | (4 – 11)(6 – 9) | (4 + 13)(6 + 15) |
| (6 + 11)(2 – 15) | (6 – 13)(2 – 9) | (6 – 9)(2 + 13) | (6 + 15)(2 + 11) |

**GoTo #5 Based on Octonion Triplet (2, 5, 7) Automorpheme: (2, 5, 7, 9, 11, 12, 14)**

| | | | |
|---|---|---|---|
| (2 + 9)(5 – 14) | (2 + 14)(5 + 9) | (2 – 11)(5 – 12) | (2 – 12)(5 + 11) |
| (5 – 14)(7 + 12) | (5 + 9)(7 – 11) | (5 – 12)(7 – 14) | (5 + 11)(7 + 9) |
| (7 + 12)(2 – 9) | (7 – 11)(2 – 14) | (7 – 14)(2 + 11) | (7 + 9)(2 + 12) |

**GoTo #6 Based on Octonion Triplet (3, 4, 7) Automorpheme: (3, 4, 7, 9, 10, 13, 14)**

| | | | |
|---|---|---|---|
| (3 + 13)(4 – 10) | (3 + 10)(4 + 13) | (3 – 14)(4 – 9) | (3 – 9)(4 + 14) |
| (4 – 10)(7 + 9) | (4 + 13)(7 – 14) | (4 – 9)(7 – 10) | (4 + 14)(7 + 13) |
| (7 + 9)(3 – 13) | (7 – 14)(3 – 10) | (7 – 10)(3 + 14) | (7 + 13)(3 + 9) |

**GoTo #7 Based on Octonion Triplet (3, 6, 5) Automorpheme: (3, 6, 5, 9, 10, 12, 15)**

| | | | |
|---|---|---|---|
| (3 + 10)(6 – 15) | (3 + 15)(6 + 10) | (3 – 9)(6 – 12) | (3 – 12)(6 + 9) |
| (6 – 15)(5 + 12) | (6 + 10)(5 – 9) | (6 – 12)(5 – 15) | (6 + 9)(5 + 10) |
| (5 + 12)(3 – 10) | (5 – 9)(3 – 15) | (5 – 15)(3 + 9) | (5 + 10)(3 + 12) |

This table was built by listing the 7 associative Octonion triplets in ascending cyclic order, then populating the first column in each such O-trip's "GoTo Listing" with the first 3 Co-Assessor pairings in its unique "triple zigzag" Trio's 6-cycle. (The cycle's last 3 pairings are identical but for internal signing.) The 1st, 2nd, and 3rd rows in the 2nd, 3rd and 4th columns are the zigzag's same-rowed (L, R) to (–R, L) Twist results. Twists based on other (trefoil) columns are (L,R) to (R, –L). The bottom row of the 2nd column's Twist is on the 3rd column's bottom, with Twists in other columns likewise on a diagonal, but oppositely oriented with respect to movements up or down the 6-cycle. The Twists for the 3rd and 4th columns, meanwhile, mirror the orientations just described for the 1st and 2nd columns, respectively. There is a "musical round" aspect to trefoils, too, one of the 3 Twists for columns 2, 3, 4 requiring its target (column 4, 2, 3) be rolled forward half a cycle. 6-cycle x 4 columns x 7 O-trips = 168 Co-Assessor couplings = PSL(2,7)'s order.

We now know that for each of the 7 O-trips there are 12 Assessors, one for each pairing of the three O-trip letters and the four O-copy Sedenions. By symmetry, the 6 pairings among these latter four units must form S-trips, and in such a manner that two uniquely match with each O-trip unit.

But each o appears in 3 automorphemes, while at the same time belonging to (28)/(7) = 4 distinct S-trips. The 6 S-trips containing the common o cannot, therefore, be distinct. But any one of 3 automorphemes sharing an o *must* share an S-trip with one of the others; moreover, it *cannot* share its *other* S-trip containing that o with that *same* automorpheme, lest they share their O-trip and be isomorphic.

Applying this to *each* pair we can pick among the three automorphemes sharing an o, leads to this

Lemma. Any two of the 3 automorphemes generating GoTo listings, and sharing an Octonion, share exactly 2 Sedenions, with the 2 each excludes from the other being shared with the *third* one. An obvious corollary is

Production Rule #3. The unique automorpheme which shares the o and S of another's Assessor will also contain that Assessor, and be the *only* distinct automorpheme to do so in the whole Sedenion space. Suppose, in the first instance, the O-trip for Assessor ( a + B ) is ( a, c, e ) and in the second, is ( a, m, n). We can solve for the Co-Assessors of the two ( a + B ) instances we know exist in the second case by setting up the usual XOR apparatus. Plug in either new O-trip unit and solve for the missing S, like this:

( a + B )( m + sg S ) ⇔ a *xor* sg S = m *xor* B; [ a *xor* m = ]n = sg S *xor* B.

Applying the other production rules to the resultant Assessor coupling yields some surprises: systems of Assessor triplicities emerge, filling the vertex structure (but only half the eight sides) of an octahedral array. These arrangements of alternating "sails" and "vents" are called ***box-kites*** in what follows. To properly appreciate them, we have one more mystery to unravel.

From our last lemma, it's clear that the 3 automorphemes P, Q, R sharing the same o collectively use up 6 of the possible 8 S's. However they are labeled, P will share 2 of its 4 S's with Q and the other 2 with R. And this is true for all automorphemes containing zero-divisors, and sharing an o: any 3 among the 7 O-trips with an o in common will therefore obey this distribution rule. In each case, there must be an S which is explicitly excluded, and XORing between it and the 3 O-trip units will generate the other 3 excludes. Is the same S excluded in all cases? Looking now at the table, the evidence is clearly yes, which leads to our

"Behind the 8-Ball Theorem". *No automorpheme which can be used to produce Assessors can contain the index-8 Sedenion.*

The table shows the "8-Ball" theorem as clearly as a proof: for GoTo #1, the O-trip is (1,2,3), and the exclude list is {8, 8+1, 8+2, 8+3}; for GoTo #2, the O-trip is (1,4,5), and the exclude list is { 8, 8+1, 8+4, 8+5 } . . . and so on. The reason for this mysterious behavior is very simple. The lower three bits of the 8-Ball are 000, making it the identity of the group of XOR operations among the Octonions' indices – the role exclusively enacted by the Real unit among the Octonions in their own domain. Hence, 8 *xor* o, o any Octonion index, will *always* be 8 + o. As any o belongs to 3 O-trips, whose units include *all* the Octonions, clearly no XORing with any other Sedenion could hope to satisfy the unit-exclusion requirements for all 7 O-trips. Indeed, for any Octonion o, none but 8 and 8 + o could do so for all three *O-trips* sharing it! Yet this is the *minimum* required, given the share-and-exclude symmetry discussed above.

The collection of seven 4-D spheres, or "balls," consisting of Quaternion copies based on the S-trip made by any o and the index-8 Sedenion, is the only such set not riven by Zero-Dividing "***slash planes***." (Not just *lines*, since *any* real coefficients attached to Zero-Dividing axes yield points which Zero-Divide each other.) This concretizes a result of Moreno's (Zero-Dividers span $2^N$-4 dimensions in $2^N$–dimensional imaginary realms, N > 3) and suggests a *unit*-based table to make "box-kites" easy to derive.

The slash-free image of the 8-Ball further suggests the Egyptian myth of Osiris, who was hacked into 14 "units" by his evil brother Set, leaving his wife and sister Isis with the Humpty-Dumpty work. The next table, then, is called the ***Osiris Partition***, and shows the labels for the 7 Octonion units running down the left, with those for the 7 viable Sedenions running across the top. The main diagonal is all null. The other 42 cells each show all four appearances, in Co-Assessor trios, of the Assessor whose indices label the cell. These are indicated by a pair of GoTo Numbers, each linked to the two Co-Assessors' indices.

The process of converting data from the Osiris Partition into box-kite drawings is a lot like a "Knight's Tour" problem. Here, the game board is 7x7 instead of the usual 8x8 used in chess, but the logic is similar. Pick a square, and write the cell's o and S numbers in parentheses as vertex labels on a triangle. In the triangle's center, put the GoTo Number. The pair of Assessor coordinates, listed in the cell with the GoTo

Number, go on the triangle's other corners. Attach another triangle, bow-tie fashion, to the first, put the second GoTo Number inside, and label its free vertices with the other Assessors' indices. Jump to the cell of an Assessor whose indices you've put on a corner, and repeat the process until you complete a "tour." You won't hit all squares, but you *will* get an octahedral vertex figure (in one of seven ways).

We will take a close look at these "box kites" in the next section. Before we begin to examine them, the reader should consider their duals, the septet of "Dunkin' Donuts" implicit in the remarks made at the bottom of the GoTo Listings table. Each 6-cycle by 4-trio array can be made into wallpaper, with the contrasting up and down rhythms of the Twist results *vis à vis* cycling direction leading to a natural representation as a toroidal two-flow system: the dunked donut moves down, the coffee's sucked up. While 6-cycling joins Assessors as vertices on box-kites, twisting joins up edges that link Co-Assessors, forming torus maps: draw a square and its diagonals, placing four Assessors from the same GoTo List and with the same o in the middle; put the four Assessors sharing the O-trip's second o in pairs at North and South; put the last four Assessors, sharing the third o, in pairs at West and East. Each half-diagonal is a join of two edges which are each other's Twist product; each edge of the square is likewise pasted to its opposite number, making the torus. Each of the four triangles defines a Co-Assessor trio residing on a different box-kite.

Neither box-kites nor dunkin' donuts can facilitate a complete Knight's Tour; but by playing these two dual systems off each other (like veins and arteries? Right brain, left brain? Ida and Pingala? Yin and Yang?), such a tour is readily conducted – even if there *is* no "group structure" to make the job description obvious!

We can take leave of the Assessors proper and consider their ensembles with <u>Production Rule #0 ("Seinfeld Production").</u> Pick a random Octonion unit o (7 ways), then a random pure Sedenion which isn't the 8-Ball or its XOR with o (6 ways); the rest follows. Or pick a random non-index-8 Sedenion (7 ways), then any o whose index isn't the XOR of that Sedenion with 8 ( 6 ways ); yadda, yadda, yadda.

# The Osiris Partition

**Two bold flush-left numbers stacked in each cell =**
**GoTo Listing numbers. (See previous table)**

| *Assessor Indices ( o + S )* | **9** | **10** | **11** | **12** | **13** | **14** | **15** |
|---|---|---|---|---|---|---|---|
| **1** | | **2** (4,15)<br>(5,14)<br>**3** (6,13)<br>(7,12) | **2** (4,14)<br>(5,15)<br>**3** (6,12)<br>(7,13) | **1** (2,15)<br>(3,14)<br>**3** (6,11)<br>(7,10) | **1** (2,14)<br>(3,15)<br>**3** (6,10)<br>(7,11) | **1** (2,13)<br>(3,12)<br>**2** (4,11)<br>(5,10) | **1** (2,12)<br>(3,13)<br>**2** (4,10)<br>(5,11) |
| **2** | **4** (4,15)<br>(6,13)<br>**5** (5,14)<br>(7,12) | | **4** (4,13)<br>(6,15)<br>**5** (5,12)<br>(7,14) | **1** (1,15)<br>(3,13)<br>**5** (5,11)<br>(7, 9 ) | **1** (4,11)<br>(6, 9 )<br>**4** (4,11)<br>(6, 9 ) | **1** (1,13)<br>(3,15)<br>**5** (5, 9 )<br>(7,11) | **1** (1,12)<br>(3,14)<br>**4** (4, 9 )<br>(6,11) |
| **3** | **6** (4,14)<br>(7,13)<br>**7** (6,12)<br>(5,15) | **6** (4,13)<br>(7,14)<br>**7** (5,12)<br>(6,15) | | **1** (1,14)<br>(2,13)<br>**7** (5,10)<br>(6, 9 ) | **1** (1,15)<br>(2,12)<br>**6** (4,10)<br>(7, 9 ) | **1** (1,12)<br>(2,15)<br>**6** (4, 9 )<br>(7,10) | **1** (1,13)<br>(2,14)<br>**7** (5, 9 )<br>(6,10) |
| **4** | **4** (2,15)<br>(6,11)<br>**6** (3,14)<br>(7,10) | **2** (1,15)<br>(5,11)<br>**6** (3,13)<br>(7, 9 ) | **2** (1,14)<br>(5,10)<br>**4** (2,13)<br>(6, 9 ) | | **4** (2,11)<br>(6,15)<br>**6** (3,10)<br>(7,14) | **2** (1,11)<br>(5,15)<br>**6** (3, 9 )<br>(7,13) | **2** (1,10)<br>(5,14)<br>**4** (2, 9 )<br>(6,13) |
| **5** | **5** (2,14)<br>(7,11)<br>**7** (3,15)<br>(6,10) | **2** (1,14)<br>(4,11)<br>**7** (3,12)<br>(6, 9 ) | **2** (1,15)<br>(4,10)<br>**5** (2,12)<br>(7, 9 ) | **5** (2,11)<br>(7,14)<br>**7** (3,10)<br>(6,15) | | **2** (1,10)<br>(4,15)<br>**5** (2, 9 )<br>(7,12) | **2** (1,11)<br>(4,14)<br>**7** (3, 9 )<br>(6,12) |
| **6** | **4** (2,13)<br>(4,11)<br>**7** (3,12)<br>(5,10) | **3** (1,13)<br>(7,11)<br>**7** (3,15)<br>(5, 9 ) | **3** (1,12)<br>(7,10)<br>**4** (2,15)<br>(4, 9 ) | **3** (1,11)<br>(7,13)<br>**7** (3, 9 )<br>(5,15) | **3** (1,10)<br>(7,12)<br>**4** (2, 9 )<br>(4,15) | | **4** (2,11)<br>(4,13)<br>**7** (3,10)<br>(5,12) |
| **7** | **5** (2,12)<br>(5,11)<br>**6** (3,13)<br>(4,10) | **3** (1,12)<br>(6,11)<br>**6** (3,14)<br>(4, 9 ) | **3** (1,13)<br>(6,10)<br>**5** (2,14)<br>(5, 9 ) | **3** (1,10)<br>(6,13)<br>**5** (2, 9 )<br>(5,14) | **3** (1,11)<br>(6,12)<br>**6** (3, 9 )<br>(4,14) | **5** (2,11)<br>(5,12)<br>**6** (3,10)<br>(4,13) | |

**{ III } Box-Kites, Lanyards, and Strut Signatures:**

The Lie algebraist's toolkit is the historical result of assimilation of imaginary units and their progeny into the orthodoxy of Number. As this century's countless failures to deal with bare charge, raw singularities of caustics, and so forth, make clear, the adjunction of Zero made by late medieval accountants to the realm of number still hasn't been fully integrated into the general population of arithmetical entities. Which makes concrete immersion in the actual phenomena being thought about all the more critical if real progress is to be made in understanding them.

The construction and "flying" of box-kites will take us a large step past these impasses. Let's start with the image of a "mock octahedron." In one $G_2$ view, they belong to a set of 7 isomorphic objects. In another, each organizes 12 sets of entities into two same-sized classes. An Assessor sits on each of the 6 vertices, each composed of two distinct lines which only zero-divide in different contexts. The vertices are all at unit distance from the origin, and have a square frame, and so on. But, due to the "quadruple covering" needed to accommodate $\pm$( A $\pm$ B ) at each, only the eight pyramids formed by each face and the origin are truly isomorphic (via "piecewise congruence") with a 3-D octahedron. The need, while circuiting along pathways joining mutually zero-dividing points, to flip from DOS to UNIX paths requires a certain *trompe d'oeil* (analogous to the hot-swapping of hard drives in a multi-protocol machine) if we're to keep the illusion going. For when these pieces are spliced together in 12-D, vertices on opposite ends of any of the three axes are a mere edge-distance apart – an Escher-like impossibility in 3-D!

The box-kite image (and octahedral box-kites are hardly unlikely) actually conveys some deep information, albeit by analogy. Call the four triangles which come from our O-trip tables "***sails***," and the four situated between them so as to prevent any two sails from sharing an edge "***vents***." In everyday spacetime, sails are 2-D sheets of silk or mylar, say, while vents allow passage of the 4-D object called ***wind***. The situation here is remarkably akin. *Our* sails live in copies of the 8-D Octonions; the Assessor trios we call ***vents***, though, require the full 16-D of "Sedenion wind" to operate in (which, as with real-world kites, is double the dimensionality of the sails).

The box-kite also has structural features which correspond to new relations among Assessors. The square frames (one per each of the three axis pairings) are a new kind of arrangement that can't exist within the confines of an O-copy, yet which are as dynamically legitimate as the trios who comprise the sails and vents. Even the struts have meaning for us, since they serve to keep diagonally opposite Assessors apart (such pairs don't zero-divide), while at the same time letting them work as higher-order components (the 4-unit additive sum of Assessors at opposite corners is itself a zero-divisor . . . and the three such composites derived one from each strut form a Co-Assessor trio).

To get started, we need some new diagramming logic – something like one might find, say, in the theory of directed graphs – to streamline our picture of how Assessors interact. And we also need to conceive some new kind of entity to formally encapsulate what these interactions comprise. As our aim is to domesticate Zero, the decimal system place-holder, let's call it a "lanyard," the string that serves as a place-holder for beads, charms, and other trinkets when building a necklace or bracelet. (Differently parsed as "LAN yard," it also suggests a collection of servers tied together to form a Local Area Network.)

The lanyard should roughly do for *multiplicative* relationships devolving upon Zero, what the formal idea of a group does for *additive* ones. The orbits of a unit's powers (or those of e which mime those) likewise form closed loops of numeric thread. Lanyards should emulate these more familiar "necklaces," offering similar service in a different context. The lanyard needs a Zero, so that B + 0 = B for all "beads" B belonging to it. And it needs a notion of sequencing akin to the Knight's Tour in its lack of redundancy and completeness of circuiting. We can't say that any two beads B and C have Zero as product; but any two beads in the lanyard must be reachable without retracings by a sequence of such links, so that adjacent beads along the way *do* have Zero products. And, we want to know the necklace has a clasp: that is to say, the lanyard is *closed*. That's equivalent to saying we can get from B to C in two directions without any recrossings – that, in programmer's jargon, the beads make a *doubly-linked list*.

Said formally, for all B, A belonging to a lanyard, there exist two indexable sets of lanyard beads E and D, of m and n elements respectively, such that: 1) $B \times E_1 = 0$, $E_i \times E_{i+1} = 0$, and $E_m = A$; 2) $A \times D_1 = 0$, $D_i \times D_{i+1} = 0$, and $D_n \times B = 0$; and, 3) sets E and D have no beads in common save for those which terminate in B and A. Finally, all Zero Divisors we know of come in pairs, so a lanyard, to be *complete*, must thread both "beads" of each Assessor – and any graphical notation used to write down lanyard structures must readily distinguish switches from UNIX to DOS paths (henceforth, "u" to "d" – as in, also, "up" or "down"). Put the indices of beads in circles or parentheses joined by lines, and in the middle of these lines write a sign: "+" if the two beads share the same internal sign, "–" if these are different.

This is motivated best by show and tell. We know that Co-Assessor trios have a 6-cycle involving changes from DOS to UNIX paths (or "down" to "up," or just "d, u" for short). Obviously, the lanyard definition was framed with things like them in mind. Using the simple graphical convention just outlined, it is easy to see why trefoils and triple zigzags must be distinguished. Consider all the "lanyard diagrams" which are possible with three links. If all three links are ( + ), we get a 3-cycle, and by the proof of Production Rule #1, we won't get lanyards at all. If we have 2 ( + ) and 1 ( – ) (and order doesn't matter) we get the trefoil pattern ( A + B => C – D => E – F => A – B ) which will cause both paths of each Assessor to be included after two circuits. 1 ( + ) and 2 ( – ) would only cover half the bases ( A + B => C – D => E + F => A + B ) and doesn't occur with Zero Divisors (and hence as a lanyard) anyway. But all ( – ) yields the triple zigzag ( A + B => C – D => E + F => A – B ), which *does* occur.

Ultimately, a lanyard *theory* should explain the mysteries, already evident in this simple 3-sign case, of which patterns make lanyards and which don't, and in which context. Such a theory might get quite interesting (and prove quite useful) in the high-dimensioned situations looming in the future. For now, though, lanyard *description* will be our modest goal. And we're about to find that there are more than two lanyard types in a box-kite. The box-kite itself forms a 12-cycle lanyard – you can get from any zero-dividing bead to any other without any retracings, which is not a trivial fact (consider the Traveling Salesman problem). And there are a handful of others, too, which we'll get to momentarily.

The first order of business, though, is to draw the box-kites as lanyard diagrams. This is rather easy, since all are isomorphic. The chart and graphics which follow summarize the facts which spell this out. And it bears pointing out that each cell of the Osiris Partition maps to *only one box-kite*. This implies that box-kites partition the Assessors without overlap.

We won't prove the isomorphy of box-kites formally: it follows immediately from placing the Assessor indices in the slots as the charting indicates, and checking against the GoTo Listings. We need, then, just a lanyard diagram of the box-kite with abstract letters filling in the slots where actual index-pairs should go, and a table of values to overwrite these with, in order to generate the box-kites themselves. Also, a stripped-down Osiris Partition, showing only which box-kites get filled by which cells, is provided as a visual aid to reflection.

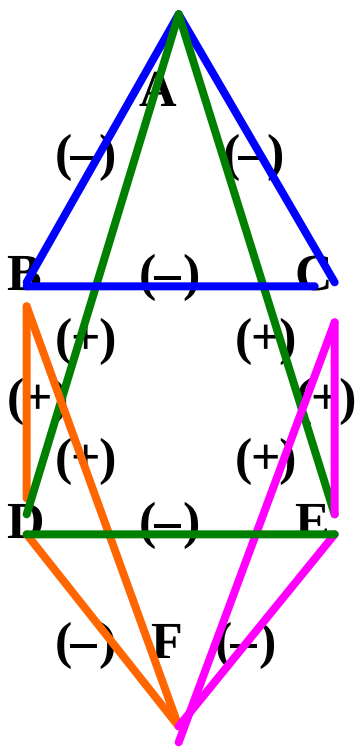

| | | | | | | |
|---|---|---|---|---|---|---|
| | III | II | V | IV | VII | VI |
| III | | I | VI | VII | IV | V |
| II | I | | VII | VI | V | IV |
| V | VI | VII | | I | II | III |
| IV | VII | VI | I | | III | II |
| VII | IV | V | II | III | | I |
| VI | V | IV | III | II | I | |

The box-kite schematic is oriented so that a triple-zigzag sail sits on top, while a triple-zigzag vent is at the bottom. Up to triangular dihedral symmetry (left-right mirroring, plus rotations in opposite senses, of A, B and C's lettering in tandem with F, E, D's) this serves to completely specify the structure. Complete consistency of letter schemes across box-kites is achieved by rotating A, B, C into O-trip cyclic order. Box-kites always have 3 trefoils and 1 triple zigzag among their sails (and also among their vents). In top-down, left-right order, the 6 Assessors of each box-kite are associated with the letters A - F in the "Strut Table" which follows, where actual index-pairs for each box-kite's Assessors are listed. (Strut pairs appear in "nested parentheses" order: AF, BE, CD.) For GoTo numbers, the triple-zigzag's is given first; the next is for the trefoil extending from the topmost Assessor to the two lowermost on either side; the two sails on the lower left, then right flanks, that join to the bottommost Assessor come last.

## STRUT TABLE

| Box-Kite | GoTo #s | A | B | C | D | E | F |
|---|---|---|---|---|---|---|---|
| **I** | 7, 6, 4, 5 | **3**, 10 | **6**, 15 | **5**, 12 | 4, 13 | 7, 14 | 2, 11 |
| **II** | 3, 2, 6, 7 | **1**, 11 | **7**, 13 | **6**, 12 | 4, 14 | 5, 15 | 3, 9 |
| **III** | 5, 4, 2, 3 | **2**, 9 | **5**, 14 | **7**, 12 | 4, 15 | 6, 13 | 1, 10 |
| **IV** | 1, 3, 5, 7 | **1**, 13 | **2**, 14 | **3**, 15 | 7, 11 | 6, 10 | 5, 9 |
| **V** | 4, 1, 6, 3 | **2**, 15 | **4**, 9 | **6**, 11 | 3, 14 | 1, 12 | 7, 10 |
| **VI** | 6, 1, 2, 5 | **3**, 13 | **4**, 10 | **7**, 9 | 1, 15 | 2, 12 | 5, 11 |
| **VII** | 2, 1, 4, 7 | **1**, 14 | **4**, 11 | **5**, 10 | 2, 13 | 3, 12 | 6, 9 |

## “Duncan Donut” for GoTo Listing #1

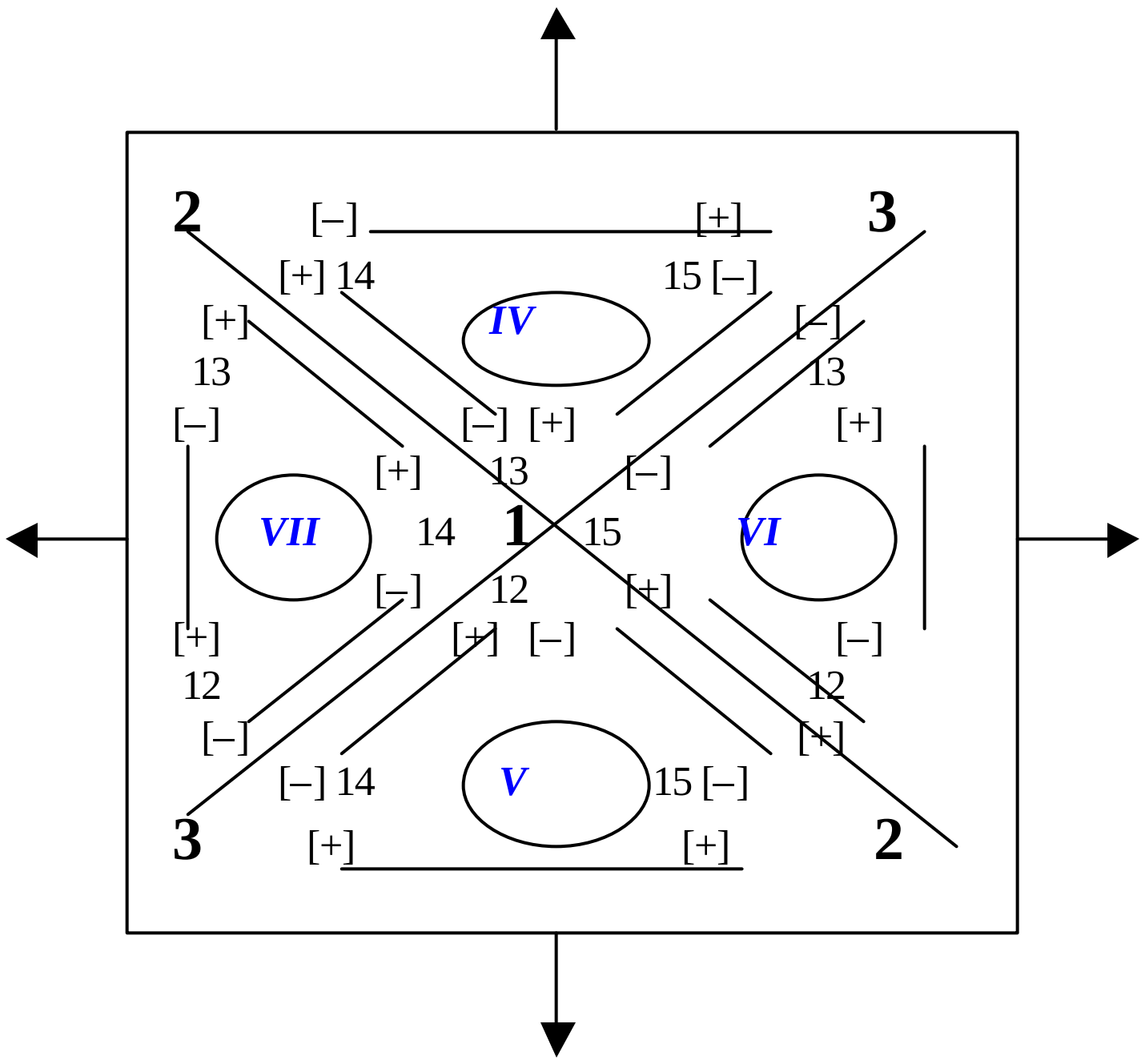

As mentioned just before the “Osiris Partition” was unveiled, the box-kites have dual structures, called “Dunkin’ Donuts” earlier. Trademark law being what it is, it’s not clear whether such a name could be deemed grounds for legal action, so the diagram above uses an alternative that sounds just the same: “Duncan Donuts,” after the great pioneer of modern dance, Isidora Duncan. Such choreography serves quite naturally as the “dual” to the “musical round” syncopations of the “GoTo Listings” which these donuts supplement. One is shown above, with arrows indicating the standard “wallpaper” method for creating a torus: join top to bottom to get a tube, then join the tube’s ends along the edges’ arrows.

The little cartouches with the Roman numerals show the indices of Trios’ box-kites. The signs in brackets tell which way to connect the two big bolded Octonion units with the indicated Sedenions along a given path. Along the bottom edge, ( 3 + 14 ) zero-divides ( 2 + 15 ), and the arrow shows they can be pasted by a Twist to the ( 2 – 14 )( 3 + 15 ) at the top. All diagonals are also borders across which Twists can be performed: ( 3 – 14 )( 1 + 12 ) in Box-Kite V is transformed by Twisting into the CoAssessor Pair ( 3 – 12 )( 1 – 14 ) in Box-Kite VII. (Readers are encouraged to draw the other donuts for themselves.)

The Roman numerals used to distinguish the box-kites are not arbitrary. They correspond to the Octonion index numbers, peculiar to each box-kite, which are called herein the ***strut signatures.*** Explaining the nature of this box-kite invariant will introduce us to the second major ingredient in Zero Divisor dynamics, after Co-Assessor trios: the paired square array of Assessor quartets, a sort of "Tray-Rack," whose opposite corners are separated by the zero-division-insulating diagonals of the ***struts***.

As is easily checked on the box-kite diagram, the Tray-Racks associated with the 3 strut pairings are 4-cycle lanyards. Each consecutive signing as we circuit the square is the reverse of those before and after, tracing a pattern suggesting the hinged collapsible frames used to support snap-on tray tables. Use the same letters with different casing, per earlier convention, for each Assessor's units, and write (b + B), (c – C), (e + E), (d – D). This is a 4-, not 8-, cycle, since *two* such *incomplete* lanyards (another starts (b – B)) must thread the same Assessors (along mutually exclusive paths) to buttress the struts. Tours around either of a Tray-Rack's paired frames never switch tracks without external interference (e.g., from a trefoil or zig-zag's impetus) – which is where the dynamic interest in their seemingly stable incompleteness resides.

Since the box-kites and their lanyard diagrams are all isomorphic, the need for deploying abstract letter variables, so invaluable in deciphering XOR signing and triplet cycles, will be much less. Instead, the next few paragraphs will frequently make use of one specific set of index values: those for Box-Kite III, which had the good fortune to be the first discovered, but is otherwise in no way special. There is much to be gained intuitively from this strategy, and (given box-kite equivalence) nothing to lose.

Putting together the generic box-kite diagram and the specifics from the last table, the Box-Kite III struts terminate in Assessor pairs indicated as { A, F }; { B, E }; { C, D }, which we can read off for the kite in question this way: { ( 7, 12 ), ( 4, 15) }; { ( 5, 14 ), ( 6, 13 ) }; { ( 2, 9 ), ( 1, 10 ) }. Each of the Assessor pairs in wavy brackets is readily shown to fail the Zero Division test. In all three instances, a single-unit product results from multiplying, and takes on one of two index values, depending on internal signing of the index pairs: 3, if signs are the same; 8, if different. Always, we're placed in the Zero-Divisor-free zone of the "8-Ball" set, with '3' being the ***strut signature*** unique to Box-Kite III.

<u>Proof.</u> (4 + 15)(7 + sg 12) = (6 + 13)(5 + sg 14) = (1 + 10)(2 + sg 9) = (1 – sg)8 + (1 + sg)3. Uniqueness is verified by visual inspection of the stripped-down Osiris Partition. The Roman numeral which names each box-kite is the only one never showing in the same-numbered row.

Note, too, the Roman numerals always appear in mirror-opposite pairs off the empty main diagonal. Each such pair in fact forms a strut, which shows the columnar layout logic of the table at the bottom: each consecutive nesting, from the outermost (A,F) pair to the innermost (C,D) is a strut, for all box kites.

Visualizations. Consider the PSL(2,7) representation of the Octonion multiplication via the logic of simplicial decomposition. Let midpoints stand in for edges, and the central point for the triangle as such. From vertex, through midpoint, to vertex, makes 3 "side triplets"; cycling through midpoints makes one "median-linking triplet"; from vertex, through center, to midpoint, makes 3 "edge-bisector triplets." Exclude the highest-dimensioned subspace (delete the central point) and the 3 edge-bisector triplets collapse, to be replaced by the binary *struts*. Joining all the non-central points yields precisely the geometry of the box-kite (with its *strut signature* now representing the index of the removed central unit). A stereographic picture of the Fano plane (PSL(2,7) *qua* simplest finite projective geometry), showing octahedral struts and four "sails" (treated as circles determined by the 3 vertices), has been published by Burkland Polster[14].

The strut signature also plays another role. Consider the "slash-product, " ( A + B)(A – B), of the two orthogonal axes which comprise an Assessor ( A, B ). For any box-kite, all slash-products of all Assessors have index 8 + o, where o is the strut signature. When we add to this the obvious fact that the square of any normed Assessor axis is the negative real unit, we get this

Lemma. All non-zero Assessor products within a box-kite are contained within the 8-Ball, and specifically within the Quaternion copy { 1, o, 8, 8+o }, where o is the strut signature. The only things of any interest in box-kite product structure, then, are the "stringing" of the Zero's in lanyards, and the stringing together of Assessors' units into higher-dimensional composites . . . leading not to Zero-Dividing orthogonal line-pairs like the Assessors, but to denser spaces we'll call *Seinfeld hyperplanes*.

* * * * *

To provide a roadmap for the remainder of this section, here's the three things left to examine before concerning ourselves with dynamics: 1) Tabulate the box-kite lanyards; 2) Collect the basic facts about Seinfeld hyperplanes; 3) Boil down our already simplistic "XOR trickery" to a degree so childish as to make higher-dimensional "box-kite" analysis easy (and give a thumbnail sketch of it for 32-D).

Let's tackle the last first and be done with it, for it's surprisingly straightforward. Consider Box-Kite VII, which has a special feature: since its strut signature is the maximal 7, the sum of the indices for each of its vertices is the maximal 15. All bits must be turned "on" after XORing to achieve this, so we needn't worry about "dropped bits." Consider the convenience this offers for higher-dimensional extrapolation. In particular, let's focus on the units in 32 dimensions, the "Path-ions" (to give them a name, after the Cabalists' "32 Paths"; in 64 dimensions, we might speak of "Iching-ions" after the Chinese system of 64 divinatory hexagrams). The one among the $2^4 - 1 = 15$ "hyper-box-kites" with maximal strut signature (of 15) would have, like all the others, $2^3 - 1 = 7$ struts. And since the sum of the indices at each vertex would be the maximal 31, we can quickly jot them down, joining strut-ends with double arrows:

| | | |
|---|---|---|
| ( 1, 30 ) ←→ ( 14, 17 ) | ( 2, 29 ) ←→ ( 13, 18 ) | ( 3, 28 ) ←→ ( 12, 19 ) |
| ( 4, 27 ) ←→ ( 11, 20 ) | ( 5, 26 ) ←→ ( 10, 21 ) | ( 6, 25 ) ←→ ( 9, 22 ) |
| | ( 7, 24 ) ←→ ( 8, 23 ) | |

Write some code to generate the 155 triplets by Cayley-Dickson process. Then, use the formula presented earlier to crank out the number of triplets for lower-dimensioned $2^N$-ions: it works for all N. (The existence and calculable number of triplets, $G_2$ derivation and 8-Ball-like generation, plus the "flexible" power law that says $X^P * X^Q = X^{[P + Q]}$ for any $2^N$-ion X, are almost all the structure we're guaranteed for indefinitely large N.) Drop anchor in the higher-dimensioned realm via the "maximal strut signature trick" just elucidated, and mysteries should begin revealing themselves in short order. Some "easy pickings": ½ (14)(13) = 91 paths between zero-divisors in each hyper-box-kite, minus the 7 struts, yields 84 Assessors unique to Path-ions in each. Add to these, the 42 inherited from Sedenions, and we have 126 – the number of roots in E7. The 84, meanwhile, cluster into 28 Co-Assessor Trios. (Recall the earlier remarks about $D_4$ and *Sedenions*, and bear in mind that in the 32-D of the Pathions, the "matter from nothing" magic of the Higgs mechanism finds its home in "Hyperdiamond Lattice" thinking.)

We know that surprises will keep on coming at least up to the $2^8$-ions, due to the 8-cyclical structure of all Clifford Algebras. And we know at least a little about what such surprises will entail: as can already be seen in low-dimensioned Clifford Algebras including non-real units which are square roots of *positive* one, numbers whose squares and even higher-order powers are 0 will appear. The 8-cycle implies an iterable, hence ever-compoundable pattern, implying in its turn cranking out of numbers which are $2^{Nth}$ roots of 0, approaching as a limit-case an analog of the "Argand diagram" whose infinitude of roots form a "loop" of some sort. If we can work with this it all, it could only be by having as backdrop some sort of geometrical environment with an infinite number of symmetries . . . suggesting the "loop" resides on some sort of *negative-curvature* surface. For the incomparably stable *soliton waves* which are deployed within such negatively curved arenas also are just about the only concrete wave-forms which meet the "infinite symmetries" (usually interpreted as "infinite number of conservation laws") requirement.

Returning to our "to-do" list, tabulating the remaining *complete* lanyards in a box-kite is easy: there are but two. First, any five Assessors form a *pyramid* atop a Tray-Rack. The pyramid's apex will have two links marked with – 's leading to "beads" bounding one of the similarly marked edges of the Tray-Rack, forming a triple zigzag; the other two "beads" of the Tray-Rack will be linked to the apex by +'s, completing three trefoils. Avoiding one of the minus-marked edges (it doesn't matter which) of the Tray-Rack, via detours through the apex, two such ¾-circuits of the Tray-Rack will yield a complete 10-cycle lanyard.

The second complete lanyard, as adumbrated earlier, is the box-kite itself. Pick one of the Tray-Racks, and call its plane horizontal. Move along the side marked *minus* which forms a trefoil with the vertical strut's upper Assessor. Mark the start and end points D1 and U2. Next, move up to the apex (the path is marked plus) and keep tracing the vertical strut's Tray-Rack edges (landing next on the node diagonally opposite U2 in the horizontal Tray-Rack) until you reach the vertical strut's lower member: U3, D4, D5. Return to the same Assessor (but opposite diagonal) you started with: mark it U6. Now, follow the *plus*-marked edge, circuiting the horizontal Tray-Rack until you reach the far end of the first leg of the journey: U7, D8, D9. Return to the upper Assessor of the vertical strut (D10), then trace out three edges of the third Tray-Rack (perpendicular to the two already traced): U11, U12, and D1 completes the cycle.

Take the two ends of the vertical strut, and join them to any edge in the horizontal Tray-Rack. Two combinations of sail-to-vent "butterflies" are possible: both triads are trefoils, or one is a triple zig-zag. In the first case, the behavior is the UUDD 4-cycle we've already seen in the Tray-Rack; in the latter case, we still get a 4-cycle, but along a different path: UUUD (or DDDU). What we have, then, in all three distinct lanyards just examined, is the predominant role of a 3 + 1 logic, itself devolving upon the 3 trefoils per triple-zigzag linkage structure which spurred thoughts about lanyards in the first place.

What this suggests, in turn, is the guiding hand, not of $G_2$, but of the Klein Group – the unsigned signature of the Quaternions, running things by remote control from the 8-Ball. And this, finally, points to another "A-D-E"-based strategy, focused on $D_4$, which underwrites the 4-D "closest packing" pattern of the 24 unit Quaternions – and also explains what might be called the "virtual geometry" of the trefoils and zigzags, whose cross-over patterns underpin box-kites. For if continuous curves are generated according to these patterns, we get diffeomorphisms of roulette curves (one from rolling an internal, the other an external, circle about another) which are the respective signatures of the Elliptic and Hyperbolic Umbilics.

These oddly named entities designate singularities in Catastrophe Theory, where A-D-E reared its head as a "Problem" in the first place. The former ("E.U.") traces a 3-cusped "control-plane" figure – a cone based on a hypocycloid where the containing circle's radius is thrice the contained – while the latter ("H.U.") traces but one cusp, generated epicycloidally. Ditto, for the lobe counts of their hypo- and epitrochoidal "behavior-plane" mappings, the former just the trefoil by another name. As Arnol'd was the first to demonstrate, these two forms of type $D_4$ are the same in 4 complex dimensions, but split into two distinct entities in 4-D real projection.

Abstractly, the Clifford algebra reading would relate the 3:1 trefoil-to-zigzag ratio to the "grading" of the line of Pascal's Triangle which spawns the Octonions – $2^3 = 1 + 3 + 3 + 1$ – with the implicit monomial terms attaching to these coefficients yielding the "interpretation" rendered above: $x^3, x^2y, xy^2, y^3$ can be read as saying "1 zigzag and 3 trefoil *sails*, 3 trefoil and 1zigzag *vents*," with x vs. y tied to the inverse "+,–" signing of strut-end links to each other Assessor. In the Catastrophist reading, these same monomials, via "Groebner bases" formalism and the symmetry of x,y interchangeability, are necessary and sufficient to indicate the writing down of germs of the unfoldings of these $D_4$ singularities: $x^2y \pm y^3$ most generally, with H.U. ("$D_{+4}$") most often written $x^3 + y^3$, and E.U. ("$D_{-4}$") usually presented as $x^3 - 3xy^2$.

These formal arguments will be given dynamic content later; for now, we still have much in the way of naturalistic description to accomplish. First off, our inventory of *complete* lanyards is done, but that of *incomplete* ones is unfinished. The full box-kite can be neatly partitioned into two halves, in two different ways, analogous to the distinct half-complete lanyards of the Tray-Rack.

First, looking at the lanyard diagram presented earlier, it's evident we can circuit its perimeter in an "oom-pah-pah" 6-cycle – a pattern, for rhythmic reasons, we might call a "Waltz Band," and wonder about the dynamic patterns of resonance between its two (UUDUUD, DDUDDU) modes. Secondly, we can navigate a circuit suggesting a "Cat's Cradle," moving only along +-signed paths. We merely go round a Tray-Rack, avoiding minus-signed edges by going up to and down from the strut orthogonal to the Tray-Rack, to get to the other +-signed side. (After which, we go *down to* and *up from* the *lower* Assessor on the strut, returning us to our starting point.) The "Cat's Cradle," too, has two modes: all U's or all D's.

And what if we attempt to do something analogous with minus-only edges? The triple-zigzags are disjoint, so there is no way to do so. What we *can* do, though, is string all 6 beads in one zigzag, and then – instead of "closing the clasp" by completing the circuit – we take a +-signed journey to the other zigzag and string all 6 of *its* beads. The result is a complete 12-bead strand with no clasp: closure in the sense a group theorist would approve would require that familiar notion from spinor theory, *double covering*. (For completeness' sake, I mention there is one other odd lanyard I'll call the "Missing Link," which covers all but one bead among the Assessors it visits: the reader is invited to discover it.)

The foremost feature of interest we have yet to investigate is Assessor composition, mandated by struts in two senses. First, precisely because Assessors at opposite ends of a strut cannot co-produce zeros, they can each contribute one or the other diagonal's pair of imaginary units into one 4-term ensemble. But this ensemble *must* zero-divide both strut-pairs of neighboring Assessors (since each zero-divides these individually), as well as *their* similarly built ensembles. *Must*, because any pair of struts necessarily span a Tray-Rack, each corner of which forms a Co-Assessor trio with its immediate neighbors. But it gets worse (or better) than this: for any strut's Assessors form trios with *all four other Assessors in their box-kite*. This implies the existence of hyperplanes of up to *four* dimensions, each point of which zero-divides each point in another of *two*. The maximal, 4-D, ones will be called "Seinfeld hyperplanes."

Use the shorthand ( A, + ) to stand for an Assessor of type ( o + S ), and (F, – ) for ( o – S ). Assume lanyard path-signing consistent with the letter scheme of the earlier box-kite diagram, and let lower-case letters stand for real scalar values. We then get these possibilities for each box-kite.

1. Lines zero-dividing strut planes: It's readily checked that composites made from strut Assessors with the same orientation zero-divide nothing, since the signs on the lanyard diagram from each to any of the other four Assessors are always opposite. Take, then, as typical of oppositely oriented generators of a strut plane this expression: m( A, + ) + n( F, – ). All such points zero-divide ( B, – ) or ( C, – ); or, ( D, + ) or ( E, + ). This gives 2 orientations x 3 struts x 4 Assessor choices = 24 ways.
2. Lines zero-dividing mismatched-edge planes: For a given Tray-Rack, bounding either +-marked edge, take as typical the mismatch p( B, + ) + q( D, – ). All such points zero-divide either ( A, – ) or ( F, + ). 2 orientations x 2 edges x 2 choices per orthogonal strut x 3 Tray-Racks = 24 ways. For either edge marked minus, take as typical the mismatch p( B, + ) + r( C, + ). All such points likewise zero-divide either ( A, – ) or ( F, + ), yielding the same count of 24 ways.
3. Strut planes zero-dividing mismatched-edge planes: Take as typical m( A, + ) + n( F, – ), which zero-divides p( B, – ) + r( C, – ) and s( D, + ) + t( E, + ). 2 orientations x 2 choices x 3 struts = 12 ways.
4. Strut planes zero-dividing strut planes: Take as typical m( A, + ) + n( F, – ), which zero-divides either p( B, – ) + t( E, + ) or r( C, – ) + s( D, + ). 2 orientations x 3 ways to pick pairs = 6 ways.
5. Triple product of three strut planes: 2 ways by orientation, but as the path joining B to C is signed plus, and that joining D to E is minus, the two struts just listed as zero-dividing the first cannot zero-divide each other without a change in orientation: we have, then, a 6-cycle trefoil- or zigzag-like relation between these three *planes*, which indicates there is in fact but *one* way to obtain this result.
6. Lines zero-dividing "tangled triad" solid: Either member from the strut in Case #1 zero-divides any picking of 3 Assessors from the 4 choices listed there. 2 orientations x 2 strut members x 3 struts x 4 ways to pick 3 from 4 = 48 ways.
7. Strut planes zero-dividing "tangled triad" solid: Use both strut members in the prior case instead of one or the other. 2 orientations x 3 struts x 4 ways to pick 3 from 4 = 24 ways.
8. Lines zero-dividing 4-space: Only a pair of strut planes can span a 4-space, leaving two candidates for the line. 2 orientations x 3 ways to pick 2 struts x 2 candidates = 12 ways.

9. <u>Strut plane zero-dividing 4-space</u>: Same choices as in Case #2, except that one of the 3 strut planes is seen as zero-dividing the 4-space spanned by the others. 2 orientations x 3 ways to pick which strut plane to leave out of the 4-space = 6 ways.

The last and ninth case subsumes all others. And as there are 7 box-kites, there are therefore 42 Seinfeld hyperplanes – the same as the number of Assessors, to which they can be uniquely mapped, in 2 ways. (The strut plane which zero-divides the Seinfeld hyperplane must have, as we saw in Case #1, one of its Assessors in "up" mode, and one in "down." The strut plane being uniquely forced, identifying the Seinfeld hyperplane with one mode or the other is tantamount to specifying the Assessor.)

If Case #9 *subsumes* all others, Case #5 – the 3-strut "triple product" – is their *generator.* This twisting of *planes* through a 6-cycle underwrites the whole "Show About Nothing" for a given box-kite. With each step along this complete lanyard, a different Seinfeld hyperplane is effectively addressed and made accessible. More than this, Production Rule #2 can be globalized in a special manner that allows each box-kite's "triple product" to communicate with each other. To see this clearly, we revert to concretizing via Box-Kite III, and consider how best to annotate these processes.

A strut which can work in a triple product has its Assessors in opposite orientations. Write one in parentheses as ( o + S ), and join it with a '+' to its partner, written ( o – S ). As a convention, since the strut signature of Box-Kite III is 3, order each strut's two o's according to their cyclic order in the O-trip beginning with 3, then order the three struts in ascending order of leading o's. The O-trips, recall, are ( 3, 1, 2 ); ( 3, 4, 7 ); ( 3, 6, 5 ). These considerations, plus lanyard signing, lead to writing the 6-cycle this way:

$$( 1 + 10 ) + ( 2 - 9 ) \Rightarrow ( 4 - 15 ) + ( 7 + 12 ) \Rightarrow ( 6 + 13 ) + ( 5 - 14 ) \Rightarrow$$
$$( 1 - 10 ) + ( 2 + 9 ) \Rightarrow ( 4 + 15 ) + ( 7 - 12 ) \Rightarrow ( 6 - 13 ) + ( 5 + 14 ) \Rightarrow$$

It's readily checked from the tables that each parenthetical term to the left of an arrow zero-divides both such terms to the right. It's also obvious that arbitrary scalar values, of arbitrary sign, could be attached to the left of each such parenthetical dyad without affecting the zero result of any such "arrow product," thereby confirming it is *planes* we are dealing with. It's clear, too, that our convention for writing down the 6-cycle makes it seem to be a triple-zigzag. However, a moment of juggling terms will convince one that other conventions of ordering terms within struts, and struts within arrows, will produce all the zeros required, while seeming to follow the signing of a trefoil – or even, of the impossible ( +, +, + ) cycle!

This clash of allowable conventions has other side-effects of interest: our "strut signature cycle-ordering" protocol, for starters, also induces the left-parenthesized o terms on either side of an arrow to have opposite orientation in their O-trip from the right-parenthesized o's. 1 x 4 ➔ + 5, but 2 x 7 ➔ – 5; for the 2nd and 5th arrows, 4 x 6 ➔ 2, but 7 x 5 ➔ – 2; for the 3rd and 6th, 6 x 1 ➔ 7, but 5 x 2 ➔ – 7.

However, this side-effect is more general than that, in two senses. First, put the strut signature "3" at the top of the standard projective PSL(2,7) triangle, and extend three rays down from it – one for each O-trip which has a "3" in it – but put them in left-to-right order in any of the six possible ways, not just in the ascending sequence of second indices per "strut signature cycle-ordering." *The same effect obtains.*

Proof. Consider the triangle as a graph of directed arrows. The "effect" can then be seen this way: the arrows between the second indices in the 3 rays flow coherently, in the sense that the arrows form a reversal-free loop, as if they were an O-trip. But the arrows of the third indices in the 3 rays also flow coherently – and *are*, in fact, an O-trip, but with reverse clockwise sense. Thus, the segment-to-segment mapping of the effect in question will always show arrows with reverse directions. This effect is clearly invariant under mirror reflection along the vertical axis of the triangle (i.e., the middle O-trip ray extending from the top). But the six possible orderings of the rays break into two sets of 3, each set of which maintains the directions of all arrows in the triangle, with one set the mirror image of the other.

But the side-effect is more general in a second sense. The pair of mirror-symmetric directed graphs just discussed are not tied to a particular strut-signature's extended rays, but emerge with similar ray extensions made from *any* of the seven indices which can sit on the nodes of the graph. They are, that is, a side-effect – a kind of signature, in fact – of the initial choice of O-trip designations and orientations from among the 480 possible (or, more precisely, of the smaller subset which require O-trip products adhere to the XOR rule). This follows from the fact that none of the points or lines in the triangle is distinguished from the rest, so moving an index to a different node without altering product structure *maintains the orientation of all arrows*. (Other indices will invariably also shift nodes, but that is irrelevant to our argument.) Such moves are just order-7 automorphisms, by a kind of "central rotation" acting transitively on all the plane's point and line sets. (Other ways of displaying the relations of the Fano plane make this more visually obvious than does the PSL(2,7) triangle representation: see Polster, cited earlier.)

The reason for belaboring the point? The Twist of Production Rule #2 can be applied twice (once to each dyad) to each strut pair. But the difference in orientation just demonstrated demands the two applications be performed in correspondingly different senses. Moreover, by showing the process concretely for Box-Kite III, we will be demonstrating it, thanks to the proven effect, for all box-kites at once. What we're about to display and justify can be thought of as either a peculiar sort of product structure, or as an algorithmic procedure for taking a "Knight's Tour" of the zero-divisor structure of Sedenion space. As the tandem applications of Production Rule #2, and the ways their results are recombined into strut pairs of all other box-kites, are highly suggestive of the ways genetic material is split and remixed during fertilization, this process is called the "recombinant divisor-navigation algorithm," or "recombinant DNA" for short.

Take two consecutive terms in the strut-product 6-cycle, writing one above the other, and according to "strut signature cycle-ordering" protocol. Twist the dyads in each of the two columns in the *same* sense – for convention's sake, let's say flip the sign of the bottom dyad's Sedenion after moving it to the top dyad, while the sign of the top dyad's Sedenion is left unchanged after relocating. Then, collect the resulting dyads *crosswise* – the left bottom dyad to the top right, the left top to the bottom right. The result is a pair of struts (the third obtainable by Production Rule #1) for the Box-Kite whose strut signature is the XOR of (either pair of) *diagonally opposite* Octonion units in the two lines being twisted.

Now, twist the diagonally opposite dyads in the two struts, in *opposite* senses – for convention's sake, let's say twist the top-left and bottom-right dyads so that the bottom dyad's Sedenion sign is flipped after moving to the top left. But twist the bottom-left and top-right dyads so that the *top* dyad's Sedenion sign is flipped after moving to the bottom left. The result is a pair of struts for the Box-Kite whose strut signature is the XOR of (either pair of) the *same-column* Octonion units in the two lines being twisted.

The three struts in a Box-Kite each have two distinct Octonion units which do not appear in the others. Pairs of these, formed one from each strut pairing, generate two distinct XOR products. When the "recombinant DNA" procedure is performed on three consecutive pairings of struts in their 6-cycle, these distinct XORings therefore produce the strut signatures *for all 6 other Box-Kites*. The group-theorist's (or chess puzzler's) dream of a complete "Knight's Tour" from any strut in any box-kite to any other (and therefore from any Zero Divisor to any other) is thereby facilitated. Now let's see all this concretely, operating on the consecutive strut-pairings of Box-Kite III.

While looking at the schematic on the next page, contemplate the way its generality for all box-kites was shown. As a proof, it is highly unsatisfying, relying on properties that are not very generic. The pattern of arrow orientations seems accidental, and is not maintained by "isomorphism" in Moreno's sense. Indeed, requiring one or more O-trips to remain fixed while we "rotate" the other lines, all without changing multiplication rules, means we're working with the 168 symmetries of PSL(2,7). But don't these merely serve to screen out all but 30 symmetrically equivalent unsigned labelings from the 7! = 5040 that are possible?

The inability to come up with a "better" proof led to the discovery of some surprising "bad news": equivalence of multiplication tables does *not* hold with Sedenions where zero-divisors are concerned. But it also points the way toward understanding the mystery of the "Moreno counterfeits" with 2 mis-signed O-trips we left hanging earlier: notions of equivalence *beneath* the multiplicative structure (and hence at the level of "graphical flow patterns") *are*, in fact, non-trivial. (Otherwise, the equivalent count of elements in PSL(2,7), henceforth "the Triangle," and the number of primitive unit zero-divisors in the Sedenions, would be "pure coincidence" – and how frequently is *that* the case in mathematics?)

*　　*　　*　　*　　*

## "Recombinant DNA" for Box-Kite III:

## An Example of the Cross-Kite Strut-Jumping Algorithm

| | | |
|---|---|---|
| Twist by column,<br>XOR o's crosswise:<br>$1\ xor\ 7 = 2\ xor\ 4 = +6$.<br>Flip signs symmetrically<br>(say, both brought to top).<br><br>( 1 + 15 ) + ( 2 – 12 )<br>( 4 + 10 ) + ( 7 – 9 )<br><br>Rearrange dyads crosswise;<br>Get struts for **Box-Kite VI**.<br><br>( 1 + 15 ) + ( 7 – 9 )<br>( 4 + 10 ) + ( 2 – 12 ) | ( 1 + 10 ) + ( 2 – 9 )<br><br>( 4 – 15 ) + ( 7 + 12 ) | Twist along diagonals,<br>XOR o's vertically:<br>$1\ xor\ 4 = -(2\ xor\ 7) = +5$.<br>Flip signs asymmetrically<br>(on same side, say left).<br><br>( 1 – 12 ) + ( 2 – 15 )<br>( 4 + 9 ) + ( 7 + 10 )<br><br>Rearrange dyads vertically;<br>Get struts for **Box-Kite V**.<br><br>( 1 – 12 ) + ( 4 + 9 )<br>( 2 – 15 ) + ( 7 + 10 ) |
| Twist by column,<br>XOR o's crosswise:<br>$4\ xor\ 5 = 7\ xor\ 6 = +1$.<br>Flip signs symmetrically<br>(say, both brought to top).<br><br>( 4 – 13 ) + ( 7 + 14 )<br>( 6 – 15 ) + ( 5 + 12 )<br><br>Rearrange dyads crosswise;<br>Get struts for **Box-Kite I**.<br><br>( 4 – 13 ) + ( 5 + 12 )<br>( 6 – 15 ) + ( 7 + 14 ) | ( 4 – 15 ) + ( 7 + 12 )<br><br>( 6 + 13 ) + ( 5 – 14 ) | Twist along diagonals,<br>XOR o's vertically:<br>$4\ xor\ 6 = -(7\ xor\ 5) = +2$.<br>Flip signs asymmetrically<br>(on same side, say left).<br><br>( 4 + 14 ) + ( 7 + 13 )<br>( 6 – 12 ) + ( 5 – 15 )<br><br>Rearrange dyads vertically;<br>Get struts for **Box-Kite II**.<br><br>( 4 + 14 ) + ( 6 – 12 )<br>( 7 + 13 ) + ( 5 – 15 ) |
| Twist by column,<br>XOR o's crosswise:<br>$6\ xor\ 2 = 5\ xor\ 1 = -4$.<br>Flip signs symmetrically<br>(say, both brought to top).<br><br>( 6 + 10 ) + ( 5 – 9 )<br>( 1 + 13 ) + ( 2 – 14 )<br><br>Rearrange dyads crosswise;<br>Get struts for **Box-Kite IV**.<br><br>( 6 + 10 ) + ( 2 – 14 )<br>( 1 + 13 ) + ( 5 – 9 ) | ( 6 + 13 ) + ( 5 – 14 )<br><br>( 1 – 10 ) + ( 2 + 9 ) | Twist along diagonals,<br>XOR o's vertically:<br>$6\ xor\ 1 = -(5\ xor\ 2) = +7$.<br>Flip signs asymmetrically<br>(on same side, say left).<br><br>( 6 – 9 ) + ( 5 – 10 )<br>( 1 + 14 ) + ( 2 + 13 )<br><br>Rearrange dyads vertically;<br>Get struts for **Box-Kite VII**.<br><br>( 6 – 9 ) + ( 1 + 14 )<br>( 5 – 10 ) + ( 2 + 13 ) |

Before becoming too attached to the prettiness of the above algorithmic pattern – providing, as it does, a very nice generalization to the "composite zero-divisor" case, of the "Seinfeld production" mode of interconnecting all *primitive* zero-divisors across box-kites -- mull the following fact. All but the simplest "lanyard theory," not to mention "Recombinant DNA" navigational strategies, vanish if we switch to an arbitrary, allegedly "equivalent" multiplication table by flipping signs on one or more units. The symmetry of what lanyard diagrams we *do* get also breaks – and breaks in different ways for different box-kites. The reader is encouraged to investigate the possibilities, and see what odd anomalies can be discovered: it is not hard, for instance, to flip a few signs and get a table which generates box-kites with at least one for which all Co-Assessor Trios are triple zigzags! (Try permuting signs on units indexed 3, 6, 5; or, try the same with 10, 15, 12 instead. Then, using the tables provided above, look at Box Kite I in either case.)

Abstractly, at least, this is not astonishing. It is known that the Cayley-Dickson process, for N > 3, generates $2^N$-ions which already suffer from a sign-flipping asymmetry of a more abstruse variety. The iterative procedure used to create double-dimensioned extensions of the reals takes one parameter (always $\pm$ 1 or 0) for each doubling, which parameter equals the square of a generating unit of the algebra. For N = 4, setting all 4 generators to the negative real unit creates the Sedenions we've been focused on. For N < 4, setting one or more generators to the positive real unit creates various modes of "split" algebras, with units some or all of whose squares (e.g., the Pauli spin matrices) can be +1, of clear utility in relativity theory. As with the closely allied Clifford algebras (which include the results of the Cayley-Dickson process), the order in which differently signed parameters are fed into the underlying algorithm makes no difference . . . *unless* N $\geq$ 4. In *that* case, as Lounesto (note 4, p. 287) points out, "in contrast to the Clifford algebras, different orderings of the parameters . . . may result in non-isomorphic algebras".

Coming closer to the ground, let's revert to the graphics lifted from Lohmus et al. earlier, and consider the underlying projective geometry of the Sedenion triplet structure. Each doubling of dimension is generated by the addition of a single point. All points are defined over the field of modulo-2 integers, so lines between points can only have coordinates whose values are 1's and 0's, meaning no line can have more than 3 points: (1,0) and (0,1) are joined by (1,1), for instance. The extra added point can hence be thought of as a "bit to the left," so that the generator of the Octonions from the Quaternions is represented by a 3-bit unit (1,0,0), while the generator of the Sedenions from the Octonions is (1,0,0,0): the "8-Ball."

By the built-in magic of projective spaces, we have dual structures: when considering projections from N dimensions to the next dimension down, the number of points is the same as the number of (N-1)-dimensional entities in the space, the number of lines equals the number of (N-2)-dimensioned things, and so on. Hence, the Octonions have 7 points and 7 lines, and Sedenions have 15 points and 15 planes. With the Octonions, there is a breakdown of algebraic symmetry manifesting in the loss of the associative property. With the Sedenions – whose 15-D continuum can be considered as spanned by units represented by the 15 points on the "projective tetrahedron" over the field of modulo-2 integers we saw before – there is a breaking of *geometric* symmetry. And it is manifested in the need to look *beneath* the usual Moreno-style "isomorphism" level, to "flow structure" patterns like those utilized in the "bad proof" above. So let's give these "not-worth-considering" properties a name, and investigate their workings.

The 168-element simple group governing the symmetries of the Triangle can be decomposed in a natural manner into 3 quite distinct subgroups, of which it is the product. The smallest – and, for immediate purposes, most significant – of these is the 4-element "Klein Group." It governs XORing among two-bit binary strings, the rotations of a tetrahedron's altitudes into each other . . . and possible label-changes on the Triangle when we keep an O-trip unpermuted and locked in place.

Let's say (1 2 3 ) is on the bottom row of the Triangle, and our "identity element" labeling scheme has 4 in the center, so that (1 7 6 ) runs up the left slope, and ( 6 5 3 ) runs down the right. Each of the 3 units in the fixed O-trip has two other O-trips running through it, and we get the other 3 elements of our Klein group by "rotating" these pairs: with '1' as pivot, the median (1 4 5 ) is swapped with the left edge, putting 7 in the middle, and 5 at the apex. Similarly, '2' as pivot interchanges the vertical and "circular" lines' other terms, swapping 6 into the center and 4 to the apex, and making 5 and 7 exchange places. The third rotation, with '3' as pivot, puts the 5 in the middle, and the left-slope midpoint, 7, at the apex.

What we're interested in doing is freely swapping any of the 7 O-trips into the bottom row where we've just put ( 1 2 3 ). This implies a 7-element cyclic group, which in turn requires the triangular dihedral group be included, to allow for rotating the Triangle and/or reflecting through median lines – which, combined with the Klein Group, has the side-effect of allowing any permutation of the "fixed" O-trip to replace it in the "fixed" position. Now let's see *why* we're interested in this.

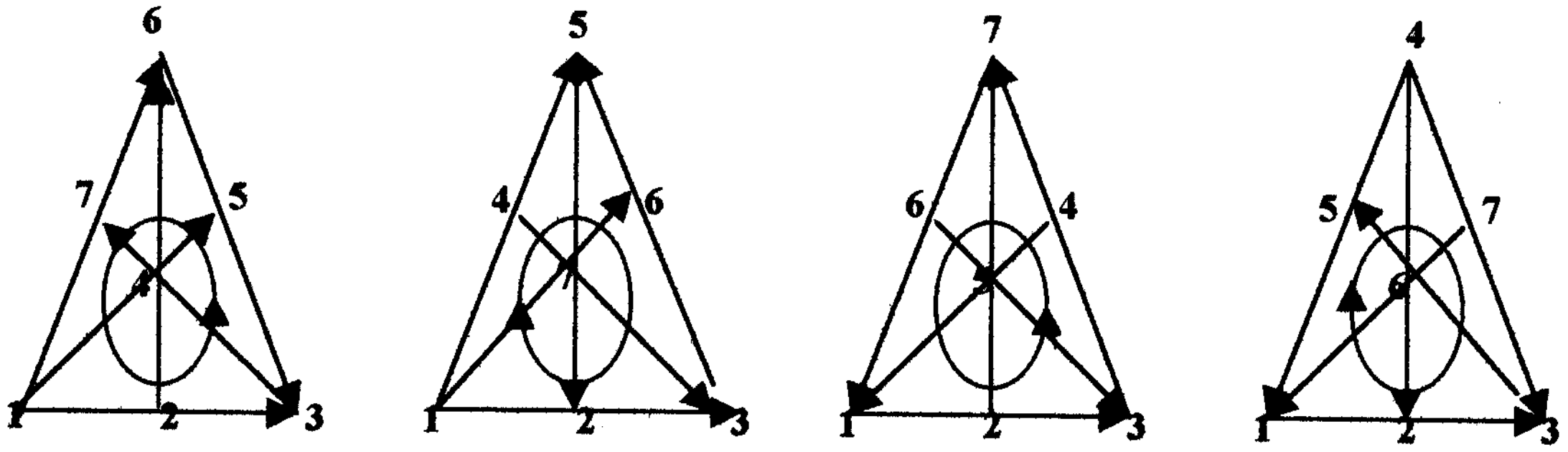

The labeling operation on the leftmost Identity Triangle which appears on the far right – the product, in turn, of the two in the middle – is the scheme of arrows implicit in the "bad proof" conducted above. Consider the "Identity" scheme, flush left: take any of the 7 O-trips, and arrange the Triangle's labels so that the chosen O-trip sits on the bottom, with the arrows of the others arrayed as shown. We say the new arrangement is "flowmorphic" to it. By using PSL(2,7)'s group operations, we can "rotate" any O-trip into the bottom-most position, so that 7 such "flowmorphs" can be produced. That's just for Octonion labelings: but recall our earlier observations concerning the "8-Ball," and you'll see that *all* labelings of the Triangle with an O-trip at bottom and '8' in the middle are "flowmorphic" to our "Identity" scheme.

For any S-trip ( o, 8, 8+o ) will have all arrows from any of the 3 o's through the '8' pointing toward the pure Sedenion. The extra bit in the '8' has, meanwhile, no effect on the orientation of the remaining 3 lines, since XORing will remove it between S's, and preserve it between an o and an S. The 7 O-trips thereby generate 7 8-Ball-containing O-copies, all "flowmorphs" of the Octonions proper.

When the mystery arose, we didn't yet know that all "O-copies" containing the 8-Ball could not harbor zero-divisors. So we see that "flowmorphs" are not merely isomorphic in Moreno's sense, but also include precisely and only those "Moreno isomorphs" which do *not* harbor zero-divisors. How to get the converse – that copies of the Triangle which *do* harbor zero-divisors are *not* flowmorphic to the Octonions proper or the septet of 8-Ball-holding Triangles in the PGL(4,2) Tetrahedron?

Let's bring to mind a built-in asymmetry of the Octonion flow-system we're focused on: if we circuit the Triangle labeled like our 8 flowmorphs, orientation along edges is not invariant, but reverses either once or twice, depending upon the clockwise sense of one's circuiting. The Sedenion tetrahedron, however, has *rotational symmetry*, which will conflict with the flow-system asymmetry of our 8 flowmorphs, yielding the "2 signing conflicts" we cited as a computable fact (and unresolved mystery) earlier.

Proof. Consider the Octonions' Triangle to be the base of the tetrahedron, then consider the face containing the (1 2 3) O-trip. But for this bottom-edge-spanning O-trip, all flow-lines on the face are S-trips included in an 8-Ball-containing O-copy. The 3 flowmorphs incorporating the O-trips joining edge midpoints to vertices in the base obviously include the 8-Ball in the center and the 12 at the top, so this covers the face's 3 lines proceeding down from the apex to the 1, 2 and 3. Similarly, the two diagonal flow-lines belong to the 2 flowmorphs incorporating the O-trips spanning the base's other two edges. Finally, the midpoint-cycle is derived from rotating the base's midpoint circle (2 5 7) onto an edge, putting the 8-ball in the middle, and getting the face's midpoint-cycle as a same-way-oriented edge.

Rotate the tetrahedron in the same sense as the (2 5 7) midpoint circle's arrow. This will cause the midpoint circles in the 3 faces to mirror the orientation of the base's ( 2 5 7 ). But it will also induce a cycling through the 3 flowmorphs straddling midpoint-to-vertex joins in the base. The flow-asymmetry of the base, however, will cause one of the 3 flow-lines shared with the face (the rightmost edge, which points down in flowmorphs) to reverse orientation, as the opposite, *up*-oriented, edge of the flowmorph is rotated through it. Meanwhile, the same sequence of rotations will cause the 3 flowmorphs set on the base's edges and with apices in *face*-edge midpoints to be cycled through, with the same result on the median descending to the '3' at bottom-right. And, since the midpoint-cycle is unaffected by the rotations, we get *only* 2 flow-lines with reversed orientation in the face. Lastly, we saw before that no flipping of signs on units can bring the face into correspondence (and hence make it flowmorphic) with the Octonions' Triangle. By symmetry, the same result obtains for all 3 faces attaching to the base, and (via rotations of the other lines of the base to its edges) to the remaining 4 8-Ball-free triangles. Q.E.D.

We can now see, too, that getting the same 2 lines to be mis-signed, no matter what inputs were fed into Moreno's mapping machine, is a red herring: ( 3 4 7 ) and ( 3 6 5 ) are just the median and edge which get reversed when ( 1 2 3 ) is the fixed O-trip in the "flowmorphic" labeling, and the Tetrahedron is rotated per above . . . and recall that Moreno's mapping specified the association of this O-trip with its first three outputs!

* * * * *

It would be easy to dwell on the last few pages' analysis and assume a negative attitude. Yet not many years ago, data later touted as evidence of Chaos Theory's ubiquity was regularly tossed out by experimentalists as proof a test had gone awry and needed redoing. In other arenas, finding – as we've just done by spinning the projective tetrahedron – that internal symmetries are injected into the behaviors you thought them insulated from, would be cause for celebration. The theory of "composed maps" in the early history of Catastrophe Theory, for instance, was built upon just such sought-for phenomena, which were sought precisely because their mathematics gave promise of mimicking the fundamental life processes of morphogenesis, where the "internal symmetries," of course, are harbored in the DNA.

In the early days of quantum theory, Niels Bohr's "Light and Life" and Erwin Schrodinger's "What Is Life?" inspired a generation of physicists to migrate to molecular biology. In our own day, some two generations after Watson and Crick, the new theory of quantum computation is leading to a rethinking about the bases of life processes that promises to inspire a second great migration. A paper delivered in the first fortnight of the new Millennium, and not yet in print (but downloadable from the Intemet) as I write this, can be taken as a harbinger of such a shift. In "Quantum Algorithms and the Genetic Code"[15], Apoorva Patel identifies components of DNA structure which implement Grover's algorithm of quantum search, proposing a key role for enzymes in maintaining the quantum coherence of the process. The view of DNA replication and protein synthesis he conjures culminates in a spectacular explanation for "why living organisms have 4 nucleotide bases and 20 amino acids": he writes, "It is amazing that these numbers arise as solutions to an optimization problem" – especially an optimization of a quantum algorithm!

What results like these and, more critically, the mind-sets which can generate them, should suggest to us is this: the usual assessment of "negative results" in studying hypercomplex numbers may require some standing on its head. Recast the problem as related to the workings of life itself, and things look somewhat different. Think of high-energy physics as cast in the image of life processes instead of the other way, and things seem more different still.

Imagine the Cayley-Dickson process as analogous to mitosis, with Sedenions representing something akin to the fourth bifurcation of a fertilized egg. The collapse of the equivalence of multiplication tables then appears as an eminently sensible defensive strategy, put in place to protect the budding embryo from external invaders: it takes on the trappings, that is, of an *immune system response*.

Without the "open sesame" of the right "keycode," isolating the lettering scheme in actual use becomes a needle-in-haystack quest. *With* the code, though, navigability via "recombinant DNA" tactics, or the earlier-discussed combination of box-kite dynamics and Twist products, allows for "making something out of nothing" in ways we've yet to fathom. What we have described herein is not likely to mark the end of the utility of "dimension doubling" in number theory; quite probably, it merely points the way past slow beginnings.

## { IV } Heuristic Afterthoughts on Dynamic Possibilities:

From the viewpoint of specialists, the realms of Lie-algebraic gauge theory and stably unfolded singularities would seem totally unrelated. Naively, however, the study of zero-divisors (arising most "naturally" in the space of Sedenions, where the Dirac Equation lives) cannot be too far removed from the study of the "zeros" of ensembles of polynomials. Philosophically, too, the only author to write major texts on both the Lie ***and*** Catastrophe theories, Robert Gilmore, finds it necessary to insert a singularity substrate beneath the smooth surface of symmetries that obsesses contemporary physicists:

> There is nothing spontaneous about symmetry breaking. It is a consequence of some dynamics or some variational principle. The word "spontaneous" is a smokescreen behind which hides a vast ignorance of the nonlinear mathematics which ultimately provides a description for the physical processes at hand. It is this nonlinear mathematics which is of central importance. The presence of symmetries allows higher catastrophes to occur at the cost of fewer control parameters. Many different symmetries may play the same restrictive role for a given nonlinear dynamics. The nonlinear mathematics is central, symmetry peripheral. Symmetry exists to be broken. To regard "spontaneous symmetry breaking" as the driving force behind some physical process is to miss entirely that physical process. ( Reference 6, p. 462 )

A triptych of figures sits immediately above the quoted passage in Gilmore's text, illustrating a sequence of phase transitions physicists often call the "Mexican hat": a parabolic bowl is poised mouth up at time zero; in the second snapshot, its bottom begins to flatten out against some plane of resistance; in the final panel, the central region where flattening has concentrated dimples up, yielding the overall contour of a sombrero. (Put metastable "false vacuum" atop its crown, "true vacuum" in the brim-rimmed trough.)

To the Catastrophist, this same set of figures signifies something apparently quite different from Higgs mechanisms governing mass's emergence from the Void. The initial picture implies a two-dimensional potential $V = x^2 + y^2$, and its squashing suggests a destabilizing perturbation, the smooth unfolding of whose effects requires a “figure of regulation” organizing the action of nonlinear “controls” on the two “behaviors” x and y. By a fairly exact analogy to the linear situation, where force is proportional to acceleration, the “frictional” effects of higher-order terms start at two derivatives’ distance from the potential in question. “Siersma’s Trick” provides a remarkably simple picture: monomials making up a regulatory ensemble or “germ” are pinpointed on Pascal’s Triangle, with those in their “shadow” at two or more derivatives’ remove comprising the “harmonics” being damped or amplified by “controls.” For the case at hand, the “germ terms” at two steps from the “behavior” at the onset are $x^4$ and $y^4$; the (nine) controls, then, reside in the 3 x 3 rhombus whose comers are the real unit at the apex (“line zero”), $x^2$ and $y^2$ at the ends of line two, and the $x^2y^2$ pendant at the center of the fourth line, of the Triangle.

From the topologist’s vantage, one of these controls is trivial (the real unit), and one is quite special (the fourth-order term). Usually called, instead, a “parameter,” this nonlinear term is special because of the peculiar invariance it underwrites. If four lines are allowed to intersect each other in an arbitrary manner (are put in “general position”) in the x,y plane, the most general ***algebraic*** representation of their situation involves the product of four terms of type ax $\pm$ by, while their most general ***geometric*** characterization would pinpoint the projective invariance of the “cross-ratio” determined by four such lines. (That is, the smooth navigation between all possible “points of view” on a scene – most concretely realized in aerial reconnaissance photography, where readily fixed signposts like railroad tracks serve to anchor the measurements – depends upon the unchanging pattern of proportionality between any four points – or, by projective-geometric “duality,” four ***lines*** – being viewed from various perspectives.)

What makes this “special” is the context in which it first emerged into view: half a century ago, Milnor and Kervair found, in 8-D and up, that manifolds could be homeomorphic yet not diffeomorphic (i.e., transformable into each other by “rubber sheet geometry,” but at the cost of loss of “smooth transitions” in their dynamic descriptions). This seemed baffling at the time, but has since become almost tame: viewed in other contexts, kindred phenomena are readily grasped – may even, in fact, seem “obvious.”

At the beginnings of perspective geometry itself, for instance, the great mystic, church leader, and mathematical visionary, Nicholas of Cusa, exploited the illusory effect of the newfangled “trompe d’oeil” painting in an inspirational discourse offered to the inhabitants of a monastery. The eyes of the painting, he said, like the “Vision of God” (his lecture’s title), seem to follow you, no matter where you stand in relation to them. If you let your ego delude you, you may come to think that God has singled you out and succumb to spiritual inflation; if you transcend this snare, however, you will get a glimpse of the higher truth that all may singly feel His gaze upon them in just the same manner, hence implying the omnivoyance of the Creator. (This illusion is now so well-known as to be a cliché, and is entertainingly discussed, bereft of all metaphysics, in a page of Steinhaus' delightful classic, "Mathematical Snapshots"; but in Cusa's day, its effect was a novelty – and a mystery – of a very high degree.)

In modern parlance, students of Chaos, Artificial Life, and systems far from equilibrium speak not of God's Omnivoyance, but of Order's "self-emergence": consider, for instance, the multitude of individual cellular beings who can combine to form the fruiting body of that singular entity we call a slime mold (an example favored by both systems theorist Ilya Prigogine[16] and catastrophist E. C. Zeeman[17]). And, instead of the "trompe d'oeil" painting, we could just hand somebody a three-mirror prismatic kaleidoscope. Unlike the ubiquitous, hexagonally-symmetric, two-mirror scopes every child knows, the three-mirror prismatic kind display an ***infinite number*** of centers of symmetry, not just one. Likewise, the Double Cusp Catastrophe (the two-potential kindred of the "sombrero" surface introduced above) can be viewed as having an ***infinite number*** of germinal centers – one for each distinct value of its "parameter."

In certain applications, this parameter can be treated as delimiting, for a small range of values at either side of its zero, a "sweet spot" within whose confines certain optimal transformations are possible which are otherwise unavailable. For "almost all" parameter values, "sweet" or not, oscillations and transformations between a very small number of stable regimes are possible (no more than five, in fact). But for the two critical values of the "parameter" which bound the "sweet spot," control breaks down utterly, and chaos (in the guise of an infinite-dimensional "control space") can ensue. This can, in fact, happen in the "real world": square sea-going vessels built with vertical sides all the way round, like floating oil-rig platforms which haven't yet been moored to the ocean's bottom, can be rendered violently unstable by small perturbations of wind and wave[18].

When $x^4 + y^4 + Kx^2y^2$ has parameter $K = \pm 2$, this ensemble can be factored into the "squared circle" $(x^2 \pm y^2)^2$. Then, instead of the figure of regulation displaying a simple polygonal symmetry (most commonly, the fourfold variety inspiring this unfolding's "Double Cusp" name), we get the limit of n-gonal symmetry as n tends toward infinity: the circular troughing, that is, of the "Mexican hat" dynamics. When, further, this general mode of unfolding is subjected to various constraints (most notably, the bilateral symmetry which controls the morphogenesis of most complex living things), the parameter can optimally take values outside the "sweet spot" – thereby guaranteeing a "Mexican hat" phase transition if the constraints are somehow suddenly relaxed, leading the untrammeled dynamic to then home in on the "sweet spot."

The possibilities just alluded to in the last paragraph can be extremely subtle, and are still not completely understood by singularity theorists. What the theorists do have a thorough grasp of, however, are the taxonomical categories: all the types of "parameter"-containing singularities (although not their detailed behaviors) have been determined up to very high dimension. And, not surprisingly, the method of their determination is a straightforward extension of the "A-D-E" logic whose cross-disciplinary ubiquity first emerged in the classificatory study of the parameter-free, "elementary" catastrophes.

Let's take the simple-minded view: for "parameter-free" catastrophes, the catalog embraces analogs of n-sided polygons (short-hand: $A_n$), their flippable siblings (the D series, "D" for "dihedral), and a finite set of three exceptional 3-D templates, the so-called Platonic solids (the E series, "E" for "exceptional" or "Euclidean," depending on taste). Meanwhile, thanks to the magic of the "A-D-E Problem," the three major division algebras can be correlated with one instance of each letter: $A_2$ for Imaginaries, $D_4$ for Quaternions, $E_8$ for Octonions. The collapse of division algebras is thereby correlated with the emergence of "parameters" and non-elementary catastrophes. Might not their "extended Dynkin diagram" logic, therefore, provide highly qualified candidates for ferreting out the secrets of zero-divisor arithmetics?

In the heart of Plato's *Timaeus*, a mathematical creation myth would spawn the building blocks of matter from the "ideal" polyhedra, but would generate these, in their turn, from three mysterious triangles. It is these which are the templates for the non-elementary catastrophes, as well – and, by "extended A-D-E" logic, for the three varieties of three-mirror prismatic kaleidoscopes also. Each "Timaeus triangle" provides the template for an infinite series of higher-order forms[19]. They can be "read off" as follows.

The singularity theorist, as if he were taking the famed analogy of "Plato's Cave" literally, determines the fundamental characteristics of his quarry in the Imaginary realm, and then projects down from complex space to our "shadow" realm of "real" solutions. As the great experimentalist, Sir David Brewster, knew soon after inventing the kaleidoscope, easy geometric arguments give only three kinds of three-mirror prismatics (and one with four mirrors which, for reasons that will be implicit in what follows, are of little interest to singularity theorists, since all the mirrors are at right angles)[20]. And, as his contemporaries were the first to fathom, ***angles*** are naturally interpreted as ***exponents*** in the complex plane.

The Timaeus Triangles, then, are readily grasped as involving three behavior variables, the component of the germ specific to each being raised to the power of the aliquot part of the semicircle indicated when putting the mirrors in place to make the prism. The equilateral triangle then leads to writing a "germ term" $X^3 + Y^3 + Z^3$, which (by the extension of "Siersma's Trick" to Pascal's Tetrahedron) has six controls (operating on X, Y, Z, and XY, YZ, and ZX) and one parameter (the XYZ monomial). The isosceles right triangle (the Double Cusp) is written thus: $X^4 + Y^4 + Z^2$ (effectively two behaviors, as the $Z^2$ part of the term can have only linear – and hence, trivial – unfolding). Finally, the 30-60-90 triangle is written as just $X^6 + Y^3 + Z^2$, with parameter of the same form as with the Double Cusp: $X^2Y^2$.

These three are respectively shorthanded as $P_8$, $X_9$, and $R_{10}$, with two less than the index indicating the number of controls. When treated as the instantiating "starter kit" for making an infinite number of higher-order forms, they are shorthanded $T_{3,3,3}$ , $T_{2,4,4}$ and $T_{2,3,6}$ in that order. We might consider the "T" to stand for "toroidal," since the polycentric nature of these singularities' unfoldings – like the invariant symmetric figure used to "wallpaper" the layout of their kaleidoscopes' ever-changing visuals – is obviously generated by a torus-map built using the appropriate Timaeus Triangle as template. (Brewster didn't use such terminology, far less that of elliptic function theory, but his pictures and geometric preliminaries convey this message quite readily to the mathematically attuned modern reader.)

Now the question of moment is how all this relates to the theory of zero-divisors in Sedenion space; the harbinger of the answer can be found in the study of "alcoves" which emerges in Tits' theory of buildings. Here, not the root structure proper, but the Weyl groups of its associated system of affine reflections, are germane. In a classic text on Coxeter theory, James E. Humphreys introduces alcoves with a telling example: if the affine Weyl group has rank 2, the associated alcoves "are triangles with angles $\pi / k$, $\pi / l$, $\pi / m$, where (k,l,m) = (3,3,3), (2,4,4), (2,3,6), in the respective cases $A_2$, $B_2$, $G_2$."[21] In other words, when moving from standard root-structure interpretation to the associated infinite-dimensional affine system, an exact correspondence is found between the simplest "non-elementary" catastrophes and the simplest non-trivial Coxeter-Dynkin diagrams – including Guillermo Moreno's favorite $G_2$!

# The *Timaeus* Tesselations (a.k.a. Brewster's Prismatics)

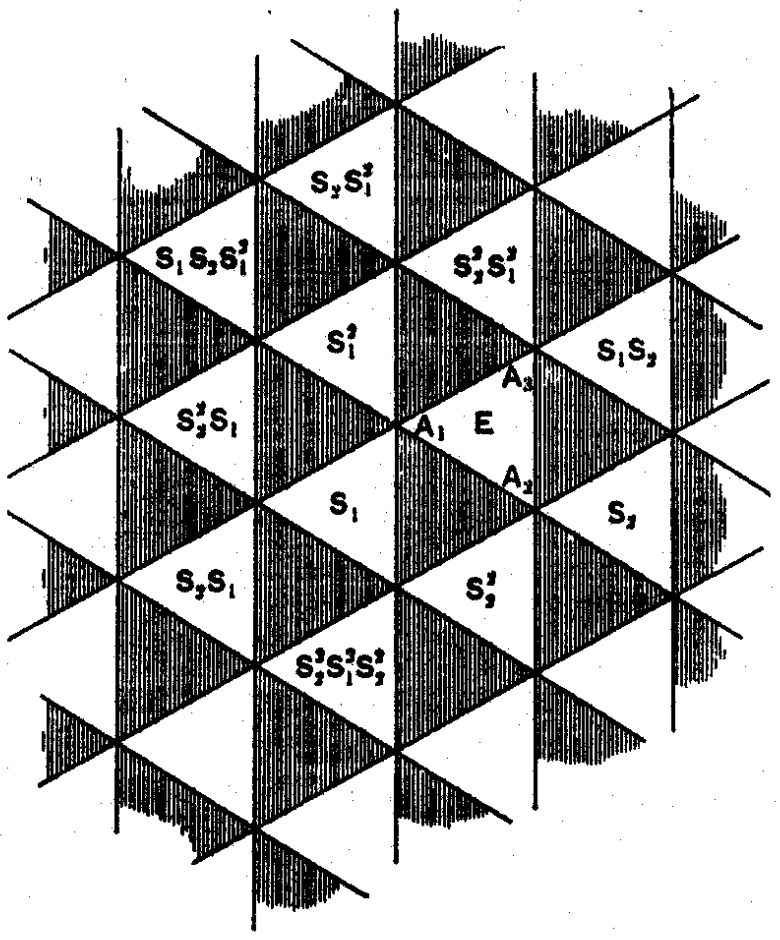

Fig. 15.

Corresponding to the solution

$$n=3, \quad m_1=m_2=m_3=3,$$

we have the general group generated by $S_1$, $S_2$, $S_3$, where

$$S_1^3=E, \quad S_2^3=E, \quad S_3^3=E,$$
$$S_1S_2S_3=E.$$

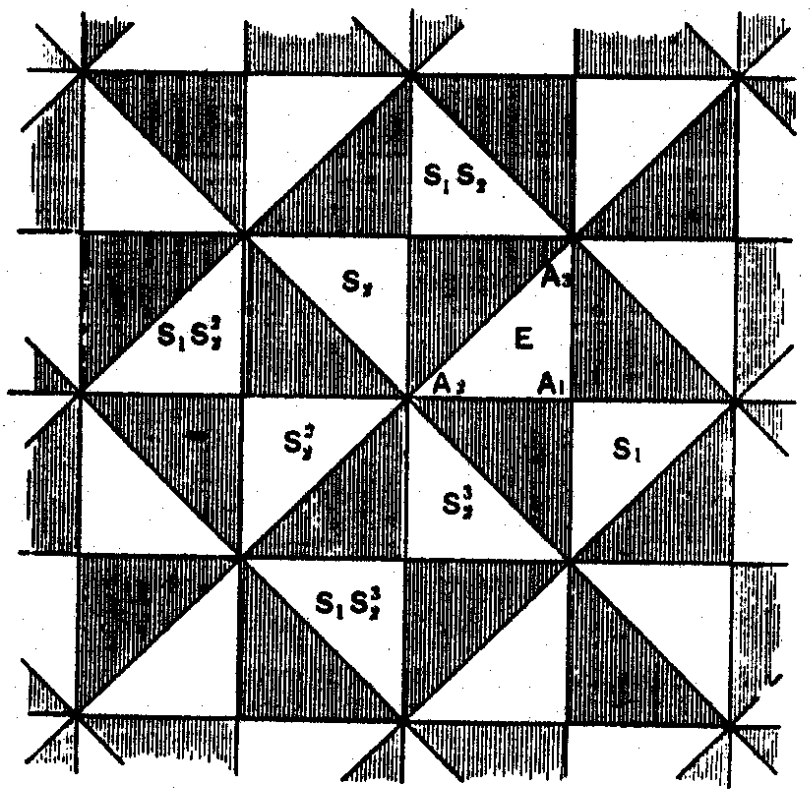

Fig. 16.

The general group corresponding to the solution

$$n=3, \quad m_1=2, \quad m_2=m_3=4,$$

is given by $S_1^2=E, \quad S_2^4=E, \quad S_3^4=E,$

$$S_1S_2S_3=E,$$

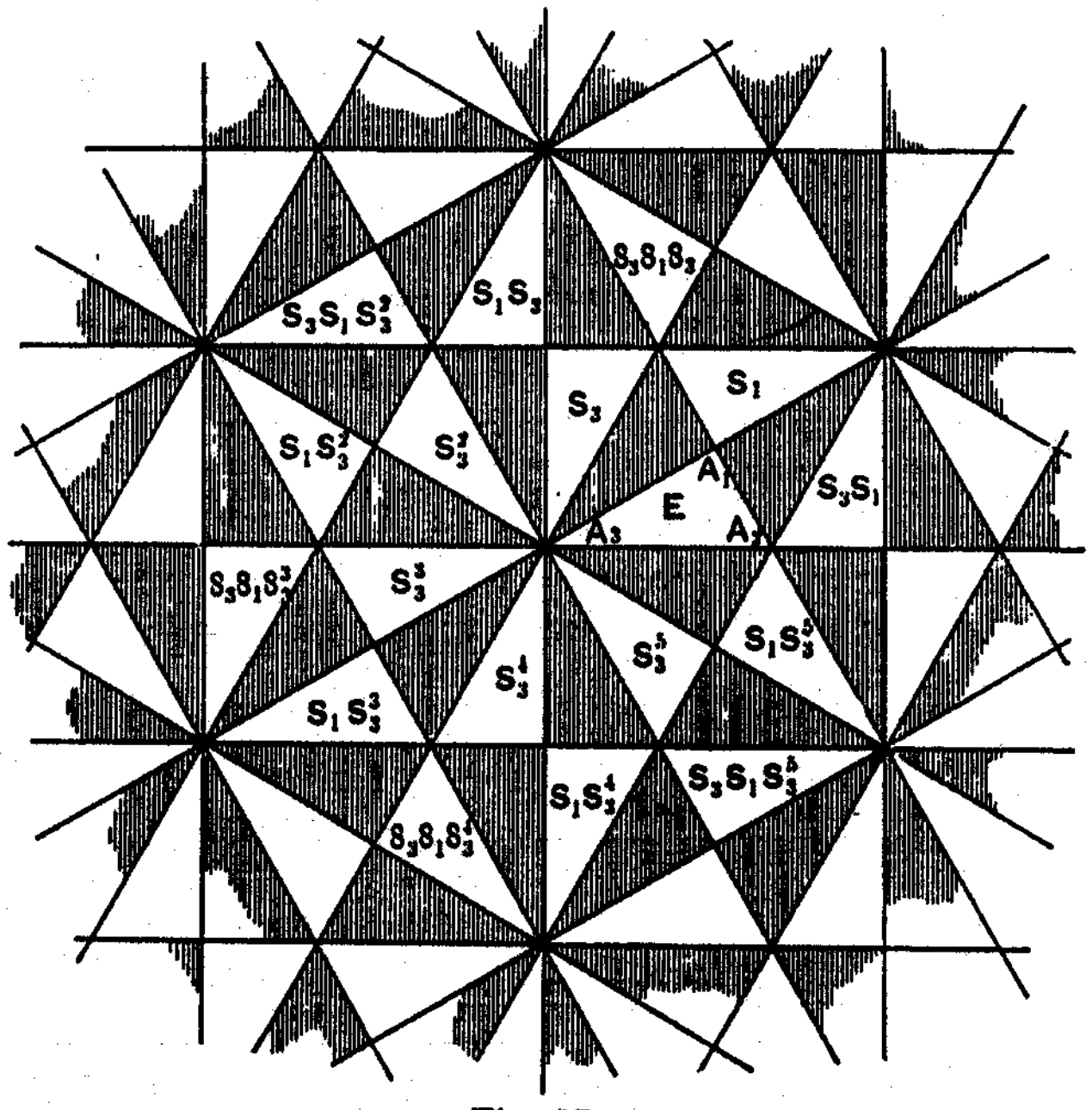

Fig. 17.

Lastly, the general group corresponding to the solution

$$n=3, \quad m_1=2, \quad m_2=3, \quad m_3=6,$$

is given by $S_1^2=E, \quad S_2^3=E, \quad S_3^6=E,$

$$S_1S_2S_3=E;$$

$G_2$ can also be construed as a singularity, either under symmetric constraint, or on the boundary of a manifold (Reference 19, pp. 278-80, where line 9 of the last page should refer to Section 17.***3***). Cast in this manner, the geometry of the germ bears a close affinity to the context of Co-Assessor Trios' toroidal "ping-pong dynamics" briefly touched upon much earlier. For now, we take some smaller steps toward a zero-divisor interpretation of Double Cusp dynamics. The reader is referred back to the graphic of the "Duncan donut" torus map, as the "accident" of its correspondence to the standard "extended Dynkin diagram" for the Double Cusp (four isosceles right triangles making a box) will be our starting point.

How could a torus map like this relate to forms of number? Considering how will take us back to the simpler (or at least, non-Sedenion) context of the 8-D Clifford algebra usually shorthanded Cl(3), and containing the Quaternions, three "Counter-Quaternion" units which are non-real square roots of ***positive*** unity ("double numbers" in the literature[22]), and an imaginary unit which commutes with all the above.

Introduced algebraically by W. K. Clifford, the Italian geometer Corrado Segré was the first to study the geometry of the number forms required by such a "double" unit. He showed its continuum of powers entails a 4-space with 1 real and 2 imaginary axes, plus one for the "double" unit itself; and, that this latter's orbit consists of two orthogonal circles centered on (mutually zero-dividing) idempotent points, passing respectively through the "double-unit" axis at ±1, and meeting in the real unit. Implicit is a kind of exponentiation defined on all points of the 2-torus "doughnut," so that its 4-space (like that of Penrose's twistor construct) is isomorphic to the product space of two standard complex circles[23].

Segré's work was long gathering dust when Charles Musès[24] independently arrived at a more modern-looking formulation, which explicitly referenced Pauli spin matrices and the operator language of the quantum-mechanical toolkit, all recast into so-called "hypernumber" forms. Musès also created an extrapolation of Octonions to 16-D, including spin-matrix analogs of the 7 imaginary units, plus the commutative imaginary now familiar in Clifford-algebra texts, but quite unusual in the physics literature of the late 60's to early 80's, when Musès wrote the bulk of his papers. This 16-D "M-algebra" has the basic structure of the superstring theorists' $E_8 \times E_8$, and so was quite prophetic of things to come a decade later. Ironically, though, the advent of theoretical approaches portended by his own led to Musès' own work falling by the wayside in mainstream research. (But Tony Smith's website has recently resurrected this rich vein of pioneering work, which the reader is encouraged to explore.)

We consider, first, the basic 4-D spin-matrix geometry. It conforms to the $X^4 + Y^4 + Z^2$ form of the Double Cusp's "Timaeus Triangle," as it: (1) is spanned by two imaginary units (period 4) and one spin unit (period 2); and, (2) has all orbits on a 2-torus. But the connection cuts deeper still. Before we can elicit it, however, we'll need a quick boil-down of Musès' notation and most basic results.

Using the Greek lower-case epsilon (**e**) for his spin units, and the usual **i** for imaginaries, he uses the subscript 0 to signify respectively the real unit and the "commutative **i**," while the subscripts 1 – 3 signify the standard Quaternion 2x2 matrices, and spin units which are equal to the same-indexed Quarternion with each of its matrix entries multiplied by the negative of the usual imaginary unit. This gives

$$e_k \cdot i_k = i_k \cdot e_k = i_0, \qquad e_k \cdot i_0 = i_0 \cdot e_k = i_k, \qquad i_k \cdot i_0 = i_0 \cdot i_k = (-e_k)$$

For $i_k$, $e_k$, $4 \leq k \leq 7$, 2x2 matrix representations are not possible, since Octonions are in general not associative; however, by constructing so-called "bi-matrices" Musès extends his e-numbers to provide a full "counter-complex" complement to the Octonions in his 16-D "M-Algebra." The indexing rules are identical to those for Octonions, so that $i_0$ commutes among same-indexed units in this context also. Hence, by the usual "associative triplets" rule combined with Musès' felicitous XORing of indices,

$$e_m \cdot i_n = i_m \cdot e_n = e_{[m\, xor\, n]} \qquad \text{and} \qquad e_m \cdot e_n = (-\, i_{[m\, xor\, n]})$$

Musès' manner of writing his "hypernumber" version of de Moivre's Theorem is rather cumbersome – in fact, he must write two separate equations, one for each circle attaching to ±e. A much cleaner and more elegant "Spin Control Equation" can be derived, relating the product of the circles to their sum:

$$(-e)^K \cdot (+e)^Q = (-e)^K + (+e)^Q - 1$$

If we set K (Q) constant, while letting Q (K) vary freely, each choice of the *constant* exponent yields a circle, isometric and parallel to the spin orbit denoted by the *varying* exponent. Most obviously, if the constant exponent be set to 0, we get one or the other spin-orbit circle; if the constant be set to ±1, we get circles parallel to the spin-orbit circles, but centered on the negative of the appropriate idempotent. For some index m, if K and Q for $(\pm e_m)$ are equal in magnitude and similarly (oppositely) signed, we get $i_0^{2K}$ ($i_m^{2K}$) – showing, again, the 4-cycle nature of the imaginaries to the spin-unit's 2-cycle. Numerous other variations can be obtained by setting the exponents to be imaginary or spin-unit quantities, etc.

Returning now to the Double Cusp, it, like all the Catastrophes, is a stratified object: it is built up entirely and only of simpler singularities, in the manner of a programmer's "object hierarchy." But as its structure is dualistic, entailing a framework for resolving conflicts between two fields of behaviors, its strata can be prepared from more than one point of view. And by containing a parameter (also called a "modulus," so that this germ is often dubbed of "unimodular" type in the literature), we can be assured its stratification has other subtleties as well. The key vantages are three:

***I: "Platonic sequence"*** -- Viewed collectively with $P_8$ and $R_{10}$ as support for the infinite T-series, $X_9$ is not only built upon a "Timaeus triangle," but contains a maximal two-behavior stratum which is associated with one of the Platonic polyhedra. $P_8$ contains $E_6$ (hence, "tetrahedrality"); $X_9$ contains $E_7$ (and thereby the symmetry common to cube and octahedron); $R_{10}$ houses $E_8$ (the dodeca- and icosa- hedral duals). Earlier, while briefly touching on the coming attractions of what we can find in the 32-D "Path-ions," it was noted that the sum total of "Assessor"-like entities in each "hyper-box-kite" was 126 – the number of roots in $E_7$. The suggestion is that this "Platonic" view of Double Cusp stratification, while pointing to ***lower***-dimensional forms in the Catastrophe context, points to ***higher*** (32-D) forms in Zero-Divisor theory. Knowing nearly nothing about the latter, exploring this possibility must be deferred.

***II: "Umbilic sequence"*** -- In discussing lanyards, the germs of the elliptic and hyperbolic umbilics ("EU" and "HU") were related to the monomials of the line of Pascal's Triangle that generated Cl(3), with their attendant geometries being related, in their turn, to the trefoil/zigzag dichotomy. In a classic study of "The Umbilic Bracelet," E. C. Zeeman[25] showed how a 3-cusped hypo- (or 1-cusped epi- ) cycloidal section could be used as the base of a cylinder which, with a 1/3 "twist," could have its two ends joined together in a torus, with this "bracelet" providing a basis for the Double Cusp's stratification. If the hypocycloid is used, all points inside the bracelet are EU's and all outside are HU's, and vice-versa for the epicycloid. The "skin," in either case, is made up of "***parabolic*** umbilic" points, whose germs unfold with one more control than EU or HU, and contain them both as strata, in precisely the manner a tilting flash-light intersects a wall in a parabola for one infinitesimal angle between its sweeping out of ellipses and hyperbolae. The closed loop traced by the cusp edges, meanwhile, supports the higher "E6" stratum.

The umbilic bracelet's dual representations are shown at left. The HU's 1-cusped epi-, and the EU's 3-cusped hypo-, cycloid sections are shown below. (From pp. 62-3, Poston & Stewart, <u>Taylor Expansions and Catastrophes</u>, Pitman, London, 1976: "symbolic" stratum = $E_6$)

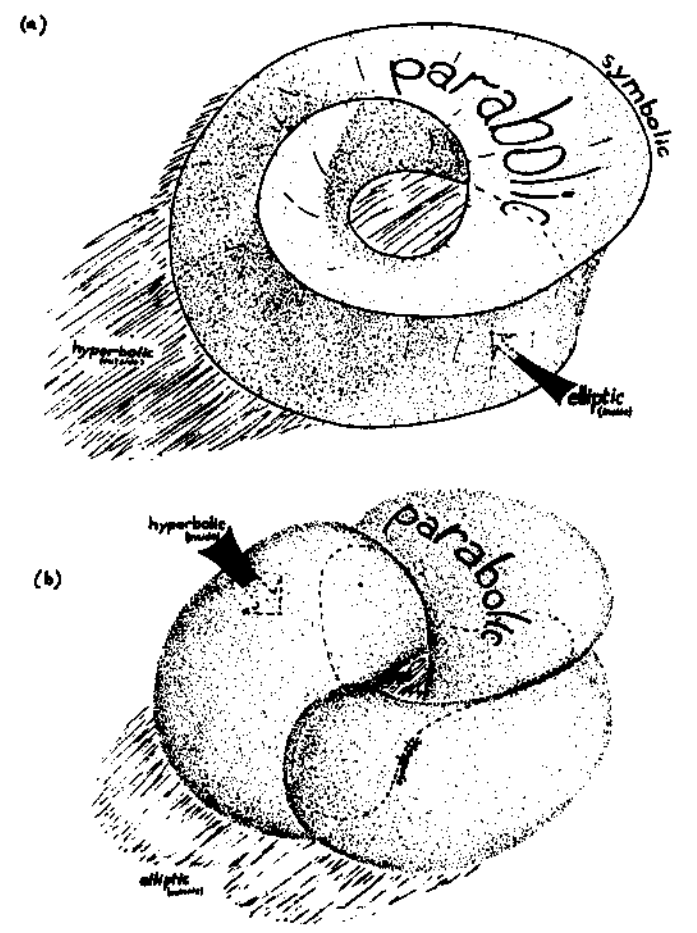

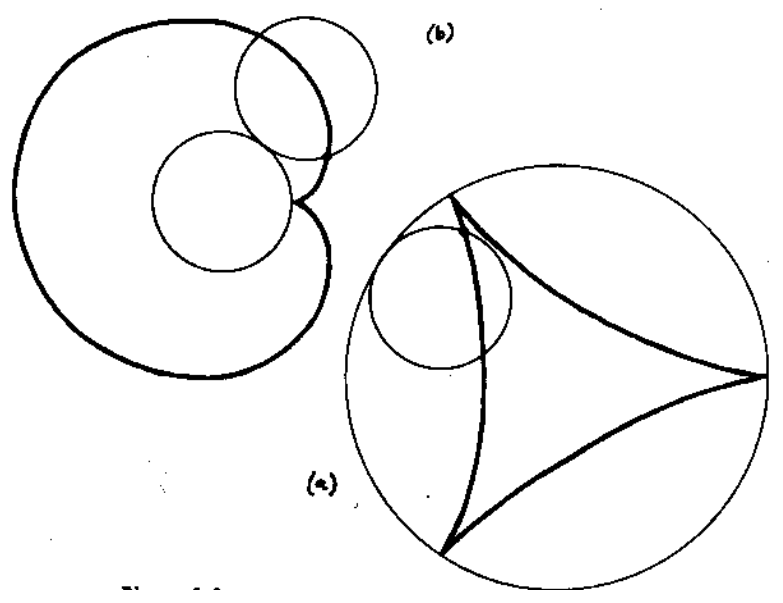

And from p. 105 of same source, "Siersma's Trick" for Double Cusp.

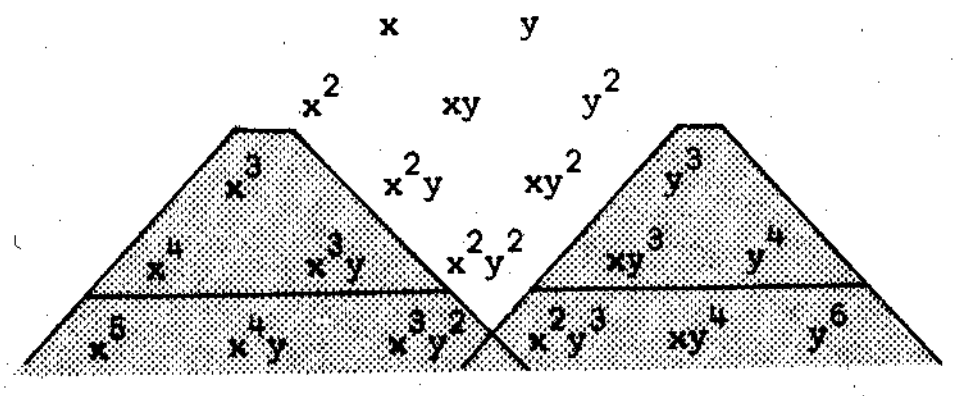

The K:Q = 2:1, 3:1, and 3:2 "power orbits" of the Spin Control Equation, displayed "six-pack" style, in unit-per-half-edge, zero-centered boxes, as described in the text on the following pages.

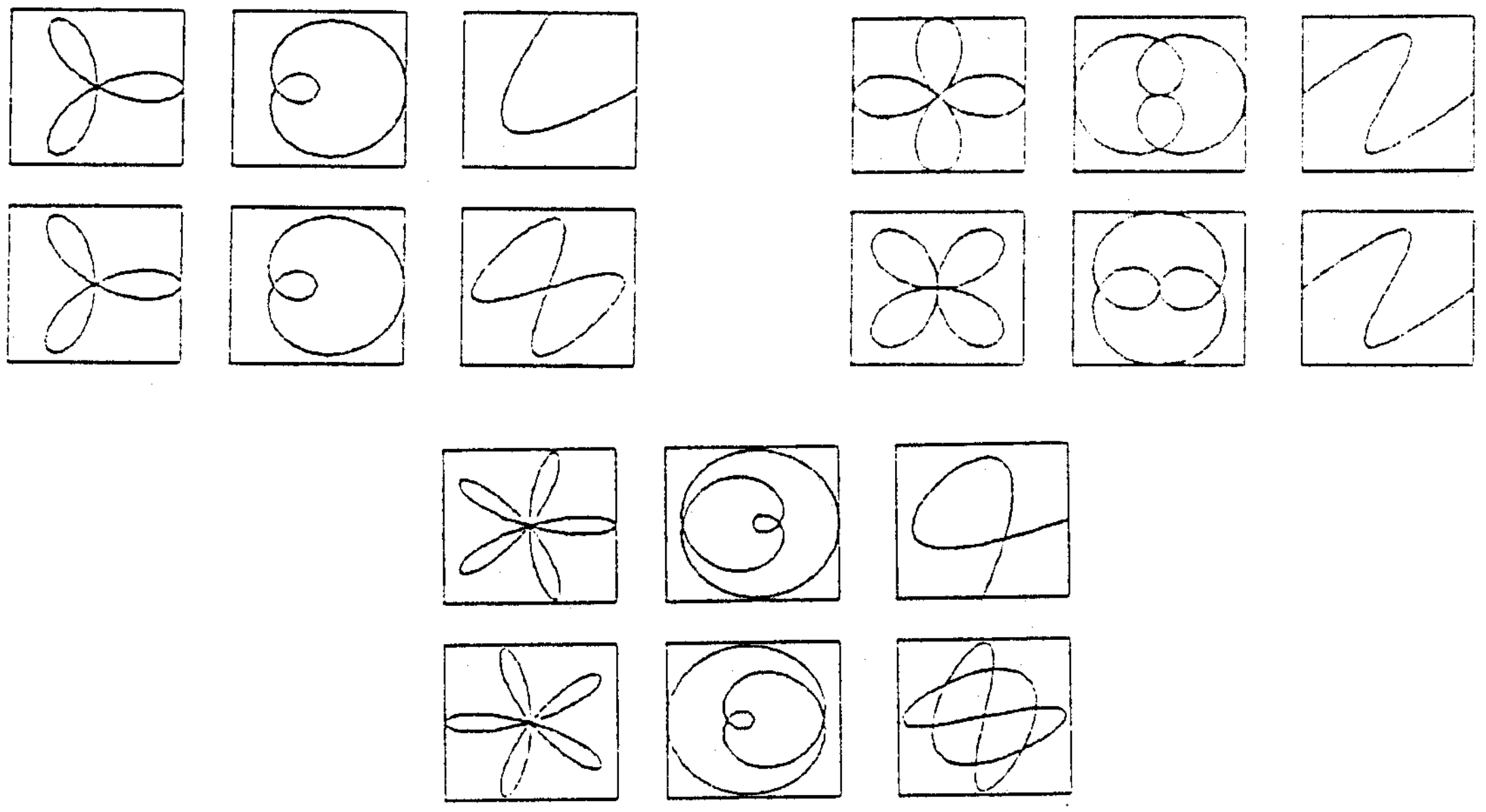

In a pair of preprints, Zeeman's associates Ian Stewart and Tim Poston built a "quartic bracelet" which tweaked out higher-order stratificational organization, devolving upon the hypocycloid of 4 cusps, or "astroid," and exploiting the cross-ratio module as basis for foliation[26]. Meanwhile, the Double Cusp's naming assumed that four minima defined the full unfolding, but James Callahan[27] demonstrated the existence of a yet higher-order stratification, allowing for 5 minima (but not necessarily a 5-cusp bracelet). This is as far as one can go, since 5 minima require separation by 4 maxima or saddles, and no unfolding of "X9" can have more than 9 critical points (which is where the form's index number derives).

This sequence emerges unbidden from the Spin Control Equation. If the two circles' exponents are in fixed rational ratio, the three simplest such proportions among relatively prime integers are 2:1, 3:1, and 3:2. (Since 1:0 or 0:1 gives us one or the other spin-unit circle, and 1:1 yields the imaginary orbits, we discount these as "trivial." We also don't care about relative signing or positioning of the exponents K and Q, as it makes no difference to our main point, and so assume both are similarly signed, and the leftmost is the larger, in what follows.) Each generates closed helical windings about the torus for its orbit; and, if we take the 4-D equivalent of an architect's plan, we get 6 planar sections, one per each pair of axes.

As a quick bit of boilerplate building with a software package like *Mathcad* readily shows, these 6 sections can always be arrayed in 3 pairs, such that each pair contains all 4 axes. The pair containing the $(1, i_n)$ and $(i_0, e_n)$ planes will show isomorphic, identically sized, rose curves, containing K + Q "petals," each housing one stable regime or "minimum" – precisely what the hypocycloids in the "control" space map into in "behavior" space. (In classical parlance dating back to Maclaurin and the first days of the calculus, the latter is generated from the former as its "pedal curve" with respect to its center.)

The pair containing the $(1, i_0)$ and $(i_n, e_n)$ planes will show isomorphic, identically sized, epitrochoids. These are curves of the broad family which Albrecht Durer first explored in the 1520's, generated by a point attached to the (possibly extended) diameter of a small circle rolling around the circumference of a larger, and corresponding, in the 2:1 case, to the map in "behavior" space of the HU's epicycloid (with the rose-to-epitrochoid pairing obtaining for the higher strata.) The remaining planes, the intersections of whose main diagonals with circles of radius $(½)^{½}$ correspond to the idempotents and their negatives -- $½(\pm 1 \pm e_n)$ -- and their square roots $½(\pm i_0 \pm i_n)$ respectively, contain Lissajous figures associated with the ratios between the spin-unit circles' angular velocities.

We know that the 2-torus in 4-D is diffeomorphic to the 3-D doughnut, so that helical "currents" are preserved in the 4- to 3-D projection. And, as we've just seen, the 6 planar sections boil down to 2 + 1 = 3 (the Lissajous figures' oscillation playing a role here akin to the dilating "third circle" used to map the 4-D orbits of the Quaternions into 3-space): i.e., to a 3-D pattern diffeomorphic to the 2-torus in 4-D.

From what was said earlier about Peixoto's Theorem, stability of toroidal currents requires there be *two* helices, with equal and rational winding numbers and opposite orientations. Assuming claims are true, we can imagine an overarching regulation of movement through these three strata having an overall symmetry of 2 x 3 x 4 x 5 = 120: the full symmetry group of the icosa- or dodeca- hedron. (Effect closure by assuming that each pair in the K:Q ratio sequence cannot be "trivial," must have its sum one greater than its predecessor, and that the larger integer in each equals or exceeds by 1 the larger integer in the prior pair. Then we get 2:1, 3:1, 3:2, 4:2 ? 2:1; that is, relative primality makes the sequence loop.)

One of the greatest surprises in singularity theory after the "A-D-E Problem" was formulated was the incredible degree to which correspondences between fields could be traced: in particular, the icosahedral reflection group just alluded to was found to govern unfoldings in the 4-D boundary of a 16-D manifold (one per each plane of symmetry, plus the identity), as well as being connected with the generic problem of "bypassing an obstacle." This latter, in turn, was shown to be anticipated – albeit not in modern language! – in the "outdated" study of the evolvents of curves pioneered by Huyghens in wave-front theory, specifically in diagrams Benniquen discovered of the evolvents of the semi-cubical parabola (itself the germ of the simplest, or Fold, Catastrophe $Y = X^3$) in the very first textbook on analysis, written by L'Hôpital from Bernoulli's lectures[28]. (As Arnol'd understates the case, "the appearance of the regular polyhedra is often unexpected.")

The former result, though, due to Lyashko[29], suggests a frontier in "zero-divisor" study that remains unexplored: what connections, if any, can be elicited between the three 16-D structures we've seen are of interest? The superstring theorist's $E_8 \times E_8$ is clearly not the same as the 16-D realization of the icosahedral reflection group's symmetry whose (symplectic) boundary is Lyashko's focus: it is just the direct product of two icosahedral rotation groups. But we know that $H_4$ (the only higher-level analogue of $H_3$) has a group structure which is the direct product of two icosahedral ***reflection*** groups, and plays a key role in the overarching theory of such five-fold singularities: Shcherbak's classic investigation of these in fact shows five-fold-symmetric forms associated with $H_4$, $H_3$, and $D_6$[30] – the last item being shown by Arnol'd to underwrite quasicrystals and their non-algorithmizable construction via the 3-D analog of Penrose tiles[31]. But this, per André Katz, lends itself to a tessellation with 22 differently decorated rhombic tiles (of either "thick" or "thin" type)[32] – corresponding to the 22 "power orbits" in Musès' M-algebra for the seven octonions, their corresponding spin-units (2 orbits each), plus one for "commutative **i**." (The third 16-D structure, of course, to interrelate with these, is that of Sedenion space itself.)

The Katz/M-Algebra construct is quite suggestive, in that "zero-divisors" are all over the place: in the space determined by $i_m \cdot e_n = e_{[m\ xor\ n]}$, m, n distinct, it's easy to see that points on the unit circle in the (m,n) plane behave peculiarly: those (as well as all others) on the diagonals through the origin are square roots of zero, and these points can, in turn, be shown to have $2^{nth}$ roots in higher dimensions (as Musès did, with unfortunately awkward notation). These diagonal pairs, meanwhile (7*6 = 42 in number, like Sedenion Assessors) are asymptotes of ***hyperbolae***, point-pairs on opposite branches of which demarcate the axes of ***elliptical*** power-orbits – a simple enough result, the search for something kindred leading to finding Yaglom's treatise on "parabolic geometry"(Reference 22.) This simple geometric support for an infinite-dimensional extension of Galilean relativity is (typically buried) undemeath the abstract language of "nilpotent sheaves" and the like in superstring theory.

The 16-D Lyashko singularity with boundary, especially given its close relationship to the "obstacle bypass" evolvent context, is quite suggestive in its own right: in particular, a novel opening to exploring some recent approaches to quantum non-localization would seem indicated. Readers are encouraged to pursue the URL and/or text version of the source containing the quotes, while keeping the "Lissajous ping-pong" motif broached much earlier clearly in mind. First, by implication, the "obstacle bypass" problem:

> One way of explaining quantum non-locality is through a hand-shaking space-time interaction between an emitter and its potential absorbers. The transactional interpretation does just this by postulating an advanced wave travelling back in time from the [future] absorber to the emitter. This interferes with the retarded wave, travelling in the usual direction from emitter to absorber to form the exchanged particle. Because both waves are zero-energy crossed phase waves, they interfere destructively outside the particle path but constructively between the emitter and absorber. The emitter sends out an offer wave and the absorber responds with a confirmation wave. Together they form a photon, just as an anti-electron (positron) travelling backwards in time is the same as an electron travelling forwards.[33]

Second, the problem of reflective boundary conditions is explicitly pinpointed as central to the same problematic. This is spelled out in the paragraph immediately following that just cited:

> Although the transactional interpretation is completely consistent with quantum mechanics, it leads to some very counter-intuitive ideas. When I see a distant quasar, in a sense the quasar radiated the photon I see long ago, only because my eye is also here to perceive it. In a sense the quasar anticipated my presence and, despite its vastly greater energy, it may not be able to radiate without the presence of myself and the other potential absorbers in its very distant future. This would imply that non-locality is observer dependent in a way which prevents any single observer having access to all the boundary conditions and hence logical prediction of the outcomes. This would prevent the universe from being computationally or deterministically predicted, but it would not prevent quantum non-locality from displaying relationships in which future states had an influence through being boundary conditions.

In summary, a third brief snippet from the same text addresses the non-algorithmicity of quasi-crystallic ordering, and its role as mediator between the chaotic and the quantized, quite succinctly:

> The quantum suppression of chaos noted in "scarring" of the wave function and other experimental results carries quantum chaos to the edge of quasi-periodicity. This is where order-from-chaos computation becomes possible, as illustrated by Conway's game of life and consistent with Freeman's model of perception.

***III: "Cuspoid sequence"*** – The simplest (because possessing but one behavior) strata are just boxcars of minima and maxima (or straddles) strung together in an A-series chain. What is surprising is just how high the strata are which the Double Cusp can contain: all of them, up through the 6-control $A_7$, which unfolds $X^8$, "are really and exactly in there, despite apparently requiring higher order terms than we have considered." (Reference 18, p. 323.) Those with even numbers of controls are "compact" in Zeeman's sense, which means that regimes can't fall off the cliff or into the sky to infinity (which, pragmatically, means the cuspoids with odd numbers of controls are often best thought of as incomplete "interrupts" requiring link-ups to something outside the ambit of the given "figure of regulation").

The two-control "Cusp" is the simplest and most famous: models of such paradoxical two-control behaviors as the rat who attacks when backed into the corner, the girl playing hard to get who succumbs only when the boy pursuing her gives up, etc., have been produced in vast numbers since Zeeman began teaching others to play this game. The maximal Cuspoid stratum can be thought of, then, as a stack of three such Cusps, whose coinciding inside one stratum may be thought of as an aspect of what a programmer would call the "UI" and not necessarily the "backend database." Put another way, the phenomena of real interest may well exist in a "three-tiered system," access to each of whose modular layers is governed minimally by a Cusp.

Each such tier-specific Cusp, meanwhile, can be thought of as having the capacity to arbitrarily bond with virtually any "behavior" beyond that defining its own dynamic, or the "umbrella dynamic" over-sighting a "handshake": which "Y" is chosen to bond with a given "X," thanks to the so-called "Splitting Lemma," is of no consequence as far as dynamics are concerned. This feature is of course absolutely critical for conceiving of software implementations allowing search and retrieval, with returns of not only precisely optimized quarry allowed, but those resulting from "bricolage" or "satisficing" as well. The Double Cusp, meanwhile, allows Chaos to intrude in certain controlled contexts, suggesting that search algorithms could employ chaos-theoretic "smart controllers," as first studied by Ott, Grebogi and Yorke in 1990, in the general case. As a pointer to interesting directions worth following, systems architect Kent Palmer has stressed the future role one can expect from "zero-divisor" number theory in the construction of any computerized modeling of complex social interactions.[34] But clearly, such issues require separate study.

"Three tiers" is the standard way of thinking about modeling sophisticated object-oriented systems these days – the three tiers corresponding to backend data, intermediary "business logic," and user interface. Each of these tiers could have its component parts scattered across any number of servers anywhere in a global network, with each tier having a life of its own, provided interfacing with the other two is built in as a given condition of the model. We've no reason to expect less than this degree of flexibility here.

An expeditious manner in which to take advantage of such potential might go like this: look for a way to "stretch" a Double Cusp (via symmetry constraints) so that the umbilic-strata machinery can provide the "message-passing" between the tiers. The starting point for such a procedure is in a classic pair of papers by Thom[35] from the early 70's, wherein the sequence of control-space figures discussed above (i.e., deltoid, astroid, pentagram) are deployed on what we'll call the "(T,V) screen" spanned by the two lowest-order controls (operating on X and Y respectively), under the constraint of bilateral symmetry.

Thom didn't treat these three as *sequentially* appearing in the same process; instead, they figure as "moments" in one 5-control equation, where one of the 3 higher controls must be "off" to effect a space-time unfolding. A logic paralleling that given earlier for collapsing the umbilic sequence to fit the contours of $H_3$ can be adopted, though, in what must be a "toy model" for now. Because of the symmetry constraint, higher-order terms than are usually accessible can be controlled; in fact, in the deltoid control-space, Thom twiddles the equations to give, in effect, an elliptic umbilic with "zoom lens": the germ of the EU ($X^3 - 3XY^2$) becomes a control, which we'll call "B," working alongside the usual "signature" control for the EU, which we'll call "C," operating on ($X^2 + Y^2$). Call this scaling or focusing effect the "offering."

In the astroid control-space, the germ of the Double Cusp itself is a control – but it is not the standard germ; rather, it includes a value of the parameter that is *outside the "sweet spot"* discussed earlier – one equal to the Pascal coefficient attaching to its monomial: $X^4 - 6X^2Y^2 + Y^4$. Call this highest-order expression's control "A." To revert to the sweet spot, the "6" must be tweaked to fall between ±2. Assume the tweaking is effected by a resonance with a Double Cusp on an adjacent tier in our three-tiered setup: an "octave," say, whose behaviors X*, Y* are such that $X^* \sim X^2$ and $Y^* \sim Y^2$, so that, when "A" is on, we get the astroid in the control plane of the two lowest-order controls (those operating on X and Y) by assuming the standard EU "C" control is "on," too. (That is, $X^2 + Y^2$ "resonates with" $X^{*2} + Y^{*2}$, with the mixed term with coefficient "6" representing, in effect, the "transaction cost.") Call this the "handshake."

In the pentagram control-space, the sequential turning off of progressively lower controls comes to a halt: only the two higher controls are "on." But there is a "mysterious change of sign" found to be required by computer simulations: the "zooming," once solicitous, now acts as a "counterpotential." This led to a correction in Thom's second paper – and suggests a heuristic interpretation: after the offering and handshake, the parties to the exchange, their business now concluded, walk away. Call this the "sign-off."

As this monograph nears its own signing off, a last question must be considered: in what manner can the above be construed as more than analogy or "toy modeling"? How, specifically, can we envision a box-kite "doing read/writes" between two Double Cusps? The secret resides in the nature of the resonance described: the "handshake" stage.

An early result of Double Cusp theory was the drawing of its "generic section" in behavior space: two intersecting, orthogonal ellipses, obeying the equation $X^2 \pm XY + Y^2 = 1$. Its direct relation to the "Holy Grail" sectioning of the astroid stratification (the astroid, as he reminds us, can be classically viewed as the envelope of a family of orthogonal ellipses) is nicely shown in Chillingworth[36]. The astroid can also be seen – as the "Holy Grail" graphic displays clearly – as the projection into 2-D of the critical edges of a curvilinear tetrahedron separating two orthogonally oriented, mirror-imaged bowls: a tellingly literal portrayal of just where the "$E_6$" is hiding inside the "$X_9$"!

Suppose each "bowl" (and hence, each ellipse) be seen as attaching to one of the two parties to the "handshake" resonance. How relate this to the box-kite formalism? The simplest route to answering this could involve another result of Charles Musès, one for which nobody, to this author's knowledge, has yet found a use. This concerns another sort of "hypernumber," which he signified by the letter "w." (He conjectured there were numerous others, which also seem to bear relation to Catastrophe Theory's visual vocabulary; the ultimate such form, he speculated, would relate to solitons: a "mad" thought which was the original inspiration for earlier remarks made herein).

Musès' "w"-numbers have power-orbits which trace out the orthogonal ellipses according to the equation in X,Y just given, but with X mapped to the reals, and Y to the "w"'s. The ellipse with a ±XY tracks the exponential powers of $(\pm w)^k$ respectively. Such orbits obey a 6-cycle (like our trefoils and triple-zigzags); and, they are, in general, non-distributive (a trademark feature of "zero-divisor" arithmetic). This much given, a method for deploying them as "hook-ups" between Double Cusps and box-kites suggests itself, using the lanyards we dubbed "tray racks." These lanyards, recall, were incomplete: they came in pairs, so that two such racks were needed to cover the "four-cornered square" such tandem structures span. Finally, note that there are three sets of tray racks, residing in orthogonal plane sections – one such section for each required mediation between the levels of our "three-tiered" software-system inspiration.

Closing question: what is the relation of the "table inequivalence" (and "immune-system responsiveness") of the Sedenions, and the "infinite-codimensional chaos at the edge of the sweet spot" of the Double Cusp? Of these two, to state vector collapse? To be pursued, perhaps, at some later date . . .

"Well, it's a long story," he said, "but the question I would like to know, is the Ultimate Question of Life, the Universe and Everything. All we know about it is that the Answer is Forty-two, which is a little aggravating."
Douglas Adams,
So Long, and Thanks for All the Fish

What is a number, that a man can know it, and a man, that he can know a number?
Warren S. McCulloch,
Embodiments of Mind

# Astroids, Holy Grails, and w-Numbers

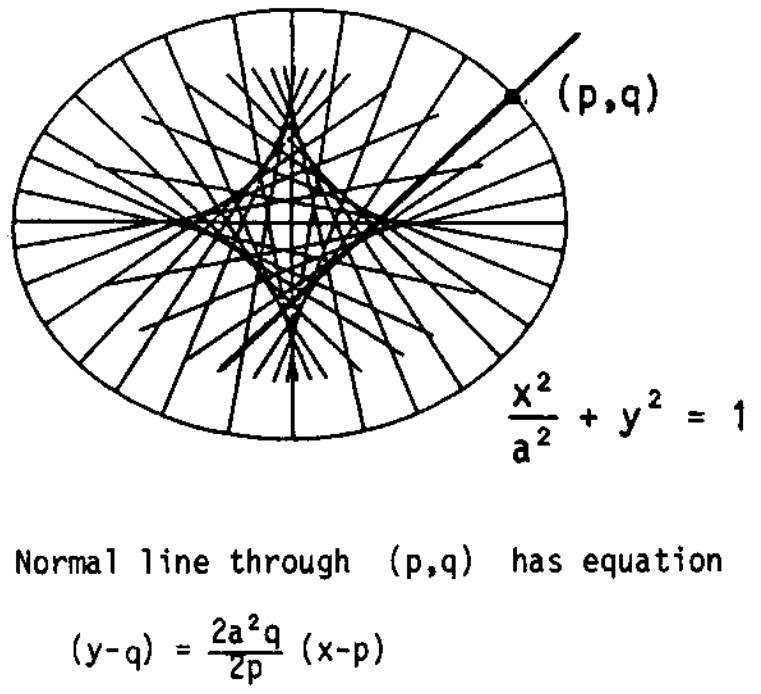

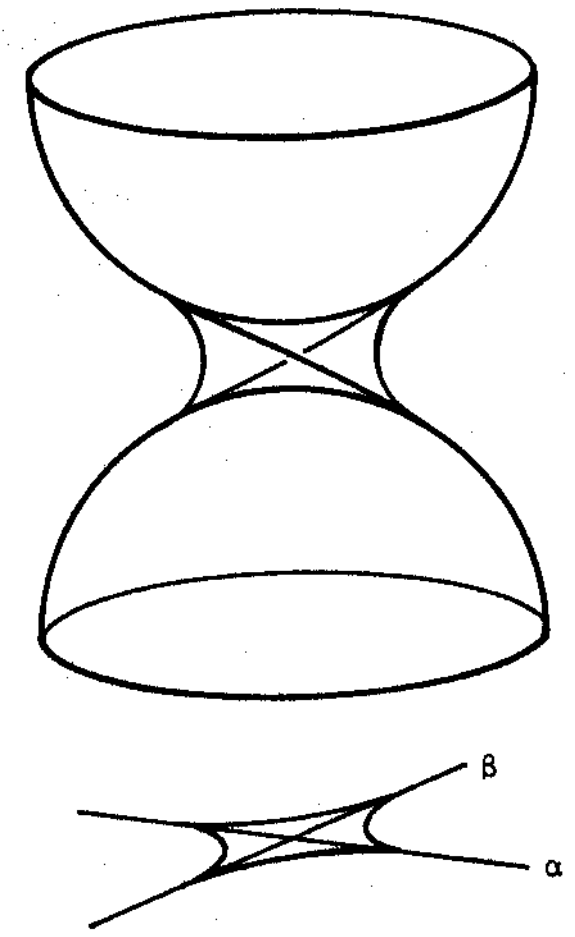

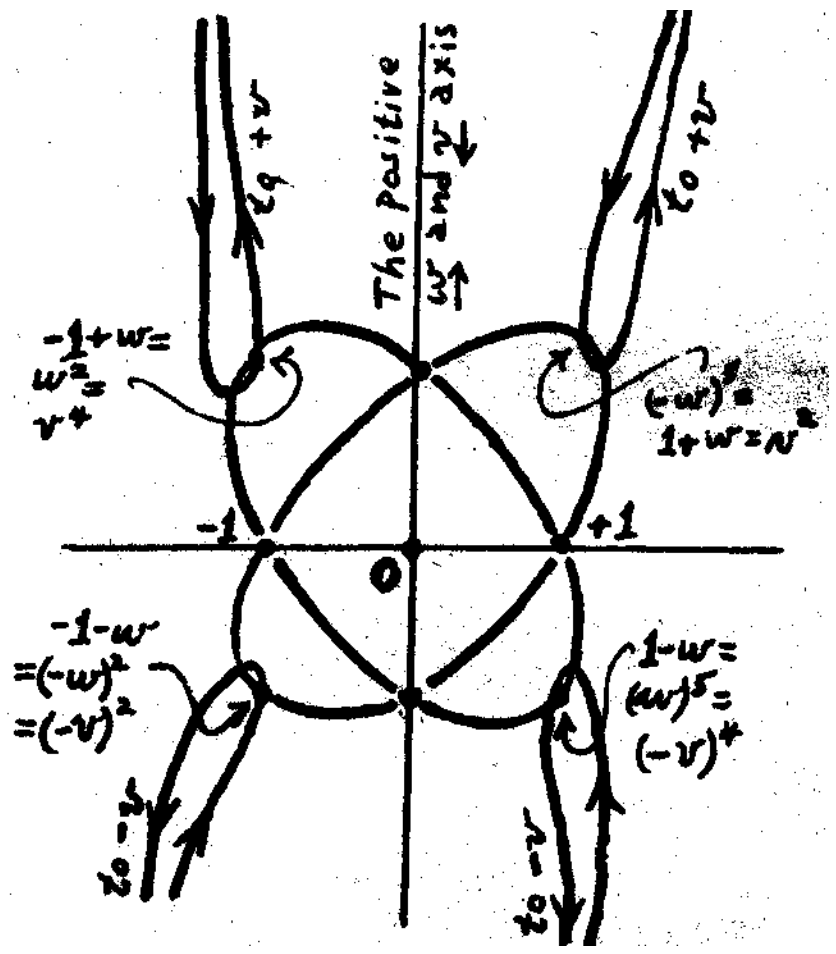

The evolute of an astroid is another, bigger astroid: the midpoints of the latter's sides are the former's cusps. The envelope of the family of normals to an ellipse is an astroid; the astroid, dually, is the envelope of a family of ellipses. In the latter case, symmetry gives us the orthogonal ellipse pair with equation $X^2 \pm XY + Y^2$: the Double Cusp's "generic section" (corresponding to a "top view" of the Holy Grail shown in the middle above, whose bowls' cross sections are our ellipses, with their tetrahedral join projecting down to yield the astroid). Meanwhile, Musès' "w-Numbers" (top right) show power-orbits of $(\pm w)^k$ traversing the same ellipses. (The infinite evolutive sequence of astroids, and $G_2$'s analogous sequence of derivation algebras for higher-order Cayley-Dixon processes, suggest a deep connection of the sort he sensed. Note, too, his "v" numbers, which just happen to attach where so-called complex "whiskers" enter the bowls of the Grail. Musès changed his mind about what these "v" might be; their only interest here is the peculiar intuition of their perpetrator: somehow he was "seeing" higher-dimensional patterns of Double Cusp strata, and sensing their relationship to patterns that cut deeper. How did he do that? And: can we run with it?)

***Offering*** ***Handshake*** ***Sign-Off***

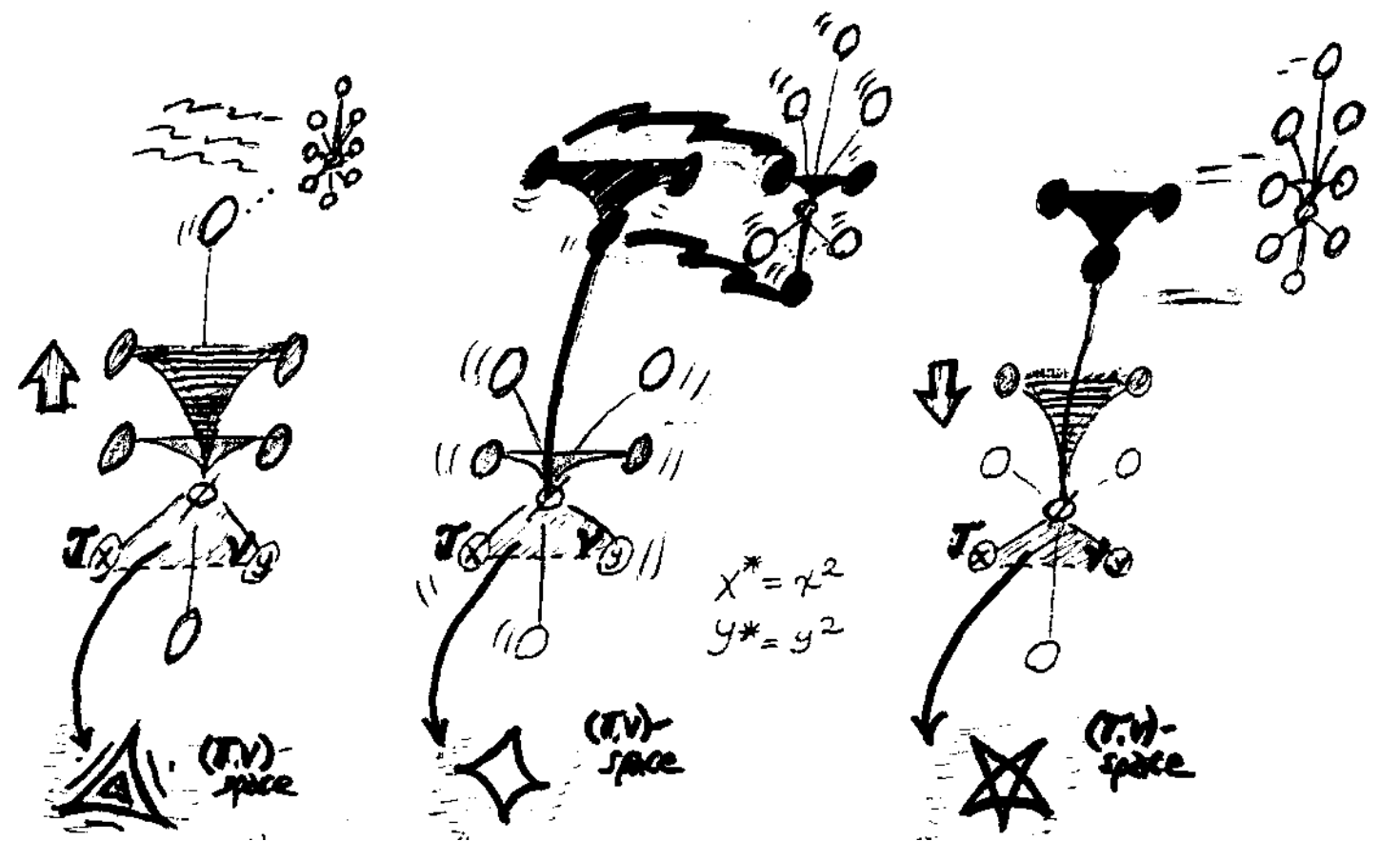

**Footnotes**

1. Guillermo Moreno, “The zero divisors of the Cayley-Dickson algebras over the real numbers,” xxx.lanl.gov/abs/q-alg/9710013, accepted for publication at Bol. Soc. Mat. Mex.
2. M. Hazewinkel, W. Hesselink, D. Siersma, F. D. Veldkamp, “The ubiquity of Coxeter-Dynkin diagrams (an introduction to the A-D-E problem), Nieuw Arch. Wisk. **25** (1977), 255-307.
3. V. I. Arnol’d, Catastrophe Theory (2nd Edition),Springer-Verlag, Berlin Heidelberg New York Toronto, 1986, esp. Chapter 16: “The Mysticism of Catastrophe Theory.”
4. Pertti Lounesto, Clifford Algebras and Spinors, Cambridge University Press, Cambridge UK, 1997.
5. E. A. Wallis Budge, The Gods of the Egyptians, vol. I, Dover, New York, 1969 (first published 1904), 418-421. One could just as well speak of the *Tibetan* Book of the Dead’s “42 Peaceful Buddhas.”
6. R. Gilmore, Lie Groups, Lie Algebras, and Some of Their Applications, John Wiley & Sons, New York London Sydney Toronto, 1974, esp. Chapter 8: “Root Spaces and Dynkin Diagrams.”
7. As of May 2000 when I write this, Tony Smith’s website URL is: www.innerx.net/~tsmith.
8. J. Lohmus, E. Paal, L. Sorgsepp, Nonassociative Algebras in Physics, Hadronic Press, 1994.
9. J. Schray, C. Manogue, “Octonionic representations of Clifford algebras and triality,” xxx.lanl.gov/hep-th/9407179. Technically, their PSL(2,7) triangle is really a representation of the isomorphic SL(3,2), which they denote PGL(3,Z2) and view as the group of the Fano projective plane over the field Z2 with 2 elements. The anomaly they pinpoint concerns the divergent behavior of “opposite” notations (e.g., all arrows relatively reversed on their PSL(2,7) diagrams) in certain contexts.
10. Susumu Okubo, Introduction to Octonion and Other Non-Associative Algebras in Physics, Cambridge U. Press, Cambridge UK, 1995. Okubo’s PSL(2,7) triangle is on p. 5.
11. W. S. Massey, “Cross Products of Vectors in Higher Dimensional Euclidean Spaces,” American Mathematical Monthly, **90** (Dec. 1983), 697-701. The original derivation of the possible cross-product rules dates to B. Eckmann, “Stetige Losunger linearer Gleichungssysteme,” Commentarii Mathematica Helvetici, **15** (1942/3), 318-339.
12. Geoffrey Dixon, “Octonions: Invariant Leech Lattice Exposed,” revised March 16, 2000, xxx.lanl.gov/hep-th/950608, and references to prior work listed therein. A good brief introduction to the Octonion cross-products can be found in Chapter 6 of F. Reese Harvey, ed., Spinors and Calibrations, Academic excesses notwithstanding, his w-numbers (and remarks in the final appendix to “Working” contrasting Eisenstein triples with the w-number 6-cycles show his intuition truly points in this direction) would seem to offer a glimpse at a profound method for codifying $G_2$ *arithmetically* – which could have enormous significance, given $G_2$’s role as the unique derivation algebra for *all* extensions of the Cayley-Dickson process.